\documentclass[dvipsnames]{article}              %  final version
\usepackage[footnotesize,bf]{caption}
\usepackage[left=1.15in,right=1.15in,top=1in]{geometry}

\usepackage[english]{babel}
\usepackage{graphicx}
\usepackage{graphics}
\usepackage{subcaption}
\usepackage{epsfig}
\usepackage{amsbsy}
\usepackage{hyperref,url}
\usepackage{amssymb}
\usepackage{amsmath}
\usepackage{mathtools}

\usepackage{chngcntr}

\usepackage{xcolor}

\usepackage{array,multirow,booktabs} % for fancier tables

\DeclareMathOperator*{\argmin}{argmin}
\newcommand{\N}{\mathbb{N}}

\newcommand{\revision}[1]{{\leavevmode\color{red}{#1}}}
\renewcommand{\revision}[1]{#1}
\newcommand{\rrevision}[1]{{\leavevmode\color{red}{#1}}}
\renewcommand{\rrevision}[1]{#1}

\newcommand{\R}{\mathbb{R}}

\begin{document}

%\volume{Volume x, Issue x, \myyear\today}
%\title[Multi-fidelity Kernel Optimization]{Kernel optimization for Low-Rank Multi-Fidelity Algorithms}
\title{Kernel optimization for Low-Rank Multi-Fidelity Algorithms}
%For at least  authors with different addresses, use instead the following commands
%\author[1]{Mani Razi}
%\author[2]{Robert M. Kirby}
%\author[1]{Akil Narayan}

\author{Mani Razi\\ Scientific Computing and Imaging (SCI) Institute, University of Utah
\and Robert M. Kirby\\ School of Computing and SCI Institute, University of Utah 
\and Akil Narayan\thanks{Corr. author: akil@sci.utah.edu}\\ Department of Mathematics and SCI Institute, University of Utah}
%\corremail{akil@sci.utah.edu}
%\corraddress{Department of Mathematics, and Scientific Computing and Imaging (SCI) Institute, University of Utah, Salt Lake City, UT, 84112}
%\address[1]{Scientific Computing and Imaging (SCI) Institute, Salt Lake City, UT, 84112}
%\address[2]{School of Computing, and Scientific Computing and Imaging (SCI) Institute, Salt Lake City, UT, 84112}
% End information for at least  authors with different addresses
% For authors with the same post address,
%\corrauthor{First A. Author}
%\corremail{f.author@affiliation.com}
%\author{Second B. Author, Jr.}
%\address{Department of Chemistry and Courant, Institute of Mathematical Sciences, New York, NY 10012, USA}
% End commands for all authors with the same address

%\keywords{kernel function optimization; low-rank approximation; multi-fidelity models; hyperparameter optimization; parameter estimation; surrogate modeling}

\date{}
\maketitle

\abstract{One of the major challenges for low-rank multi-fidelity (MF) approaches is the assumption that low-fidelity (LF) and high-fidelity (HF) models admit ``similar'' low-rank kernel representations. Low-rank MF methods have traditionally attempted to exploit low-rank representations of \emph{linear} kernels, \revision{which are kernel functions of the form $K(u,v) = v^T u$ for vectors $u$ and $v$}. However, such linear kernels may not be able to capture low-rank behavior, and they may admit LF and HF kernels that are not similar. Such a situation renders a naive approach to low-rank MF procedures ineffective. In this paper, we propose a novel approach for the selection of a near-optimal kernel function for use in low-rank MF methods. The proposed framework is a two-step strategy wherein: (1) hyperparameters of a library of kernel functions are optimized, and (2) a particular combination of the optimized kernels is selected, through either a convex mixture (Additive Kernel Approach) or through a data-driven optimization (Adaptive Kernel Approach). The two resulting methods for this generalized framework both utilize only the available inexpensive low-fidelity data and thus no evaluation of high-fidelity simulation model is needed until a kernel is chosen. These proposed approaches are tested on five non-trivial real-world problems including multi-fidelity surrogate modeling for one- and two-species molecular systems, gravitational many-body problem, associating polymer networks, plasmonic nano-particle arrays, and an incompressible flow in channels with stenosis. The results for these numerical experiments demonstrate the numerical stability efficiency of both proposed kernel function selection procedures, as well as high accuracy of their resultant predictive models for estimation of quantities of interest. Comparisons against standard linear kernel procedures also demonstrate increased accuracy of the optimized kernel approaches.
}

\section{Introduction}\label{sec:intro}
Multi-fidelity (MF) approaches have become an active area of research in the literature for uncertainty quantification and several approaches have been introduced with widespread application in many areas of science and engineering. For example, such approaches have been utilized for molecular dynamics simulations, chaotic motion, finite element analysis of acoustic problems, finite volume analysis of heat driven cavity flow, computational aerodynamics, particle-laden turbulence, multidisciplinary design optimization of wings, airfoil optimization, and plasmonics; see \cite{razi2018MD,narayan2014stochastic,zhu2014computational,hampton2018practical,skinner2017evaluation,jofre2018multi,allaire2014mathematical,lam2015multifidelity,razi2018quantized}, respectively. MF models are computational meta-algorithms that are used to simultaneously leverage the strengths of low-fidelity models (i.e., low cost) and high-fidelity models (i.e., high accuracy). 
%lowtake advantage of the potential relationship between computational models of different fidelity levels to achieve both optimal accuracy and computational efficiency.    

The setup we consider in this paper is that we are given an ensemble of models, each of which has an identical set of parameters. \revision{Note that this can be a strong assumption, as many times differing models have differing parameter sets. \rrevision{Situations when models have different parameter sets are common, but} we demonstrate that there are numerous examples in scientific computing when our assumption is valid.} Trusted but computationally expensive models are called \emph{high-fidelity} (HF) models. Relatively inexpensive computational models with low accuracy (relative to some baseline truth or HF models) are \emph{low-fidelity} models. These characteristics that give rise to low versus high-fidelity models usually reflect a trade-off between computational cost and accuracy. For example, these model ensembles are frequently produced by varying a discretization coarsening parameter, by making geometrical simplifications, or by making simplifications in physics.

Multilevel Monte Carlo-type methods are one popular class MF methods, \revision{and several other types of strategies exist} \cite{giles_multilevel_2015,fernandez-godino_review_2016,peherstorfer2018survey}. These methods are developed based on the approximation of correlation between low- and high-fidelity models with a type of reproducing kernel, or Gramian matrix. \revision{Kernel functions,  which are used to construct Gramian or covariance matrices for different applications are in essence high-dimensional generalized inner product functions.} These functions are applied to map the original the original data into a higher-dimensional spaces and often acts as the engine for automatic feature engineering. In this paper we consider a different branch of MF methods: low-rank MF methods \cite{narayan2014stochastic}. Low-rank MF methods compute an inexpensive low-rank approximation to HF models thorough exploration of the low-fidelity model. \revision{Here, the low-rank property means that emulators have a small number of degrees of freedom.} 

In the context of model discrepancy, ``fidelity'' can also be interpreted as different values of quantized parameters that produce drastic changes in the way simulation models behave. This new form of quantized fidelity level has been recently studied for problems such as different particle numbers in plasmonic nano-particle arrays, different number of species for competitive ecological systems, and different levels of density for associating polymer networks~\cite{razi2018quantized}. Although the low-fidelity models for most of aforementioned examples are typically obtained by computational simplification, data-fitted models can also be used in the context of multi-fidelity models. 

The core idea in low rank multi-fidelity methods is to utilize what is called the ``Kernel Trick'' in machine learning, which uses kernel functions to map model realizations to a kind of correlation. Until now, only linear kernels have been used in the literature for low-rank MF methods~\cite{zhu2017multi,keshavarzzadeh_convergence_2019,razi2018MD,narayan2014stochastic,razi2018quantized,zhu2014computational}, although alternative kernels are popular for other kinds of MF approaches \cite{perdikaris2015multi,perdikaris2016multifidelity}. \rrevision{Linear kernel functions are those of the form $K(u,v) = v^T u$ for vectors $u$ and $v$.} The use of the linear kernel function in low rank MF methods empirically produces accurate predictions. However, use of this kernel in MF methods has three major drawbacks:
\begin{itemize}
  \item The rank of the kernel matrix equals the size of the model realization vectors. The result is that the number of usable HF simulations is limited to this rank.
  \item The linear kernel matrix can be ill-conditioned. Although such an event could indicate exploitable low-rank structure, it also results in inaccurate finite-precision computations. This limitation can sometimes require users to employ \textit{ad hoc} approaches \cite{razi2018quantized} to mitigate the stability issue.
  \item Low-rank structure of models may not be detectable using a linear kernel, but is perhaps detectable using an alternative kernel.
\end{itemize}
In this paper we propose and empirically explore a new procedure for selecting a kernel in low-rank MF methods. 
Kernel function selection strategies have been already used in the context of predictive estimation of the models' fidelity, e.g., \cite{mehmani2018concurrent}. 
To the best of the authors' knowledge, the details of the approach employed in this paper are distinct from existing kernel learning/selection strategies in the literature. The procedure in this paper is motivated by low-rank MF algorithms, which have different desiderata than standard kernel learning approaches. However, our approach does exhibit similarities to existing methods: we optimize hyperparameters in kernels, and optimize over multiple families of well-known kernels. Overall, our approach performs this optimization by attempting to minimize the discrepancy with respect to the linear kernel matrix, and a numerical rank regularization is employed to promote stability. We propose two approaches: The \emph{Additive Kernel Approach} constructs a kernel via a convex combination hyperparameter-optimized kernel functions. The second, the \textit{Adaptive Kernel Approach} selects a single hyperparameter-optimized kernel from the library using an inexpensive high-fidelity verification. Both new approaches perform optimization on LF data, and do not require any additional HF evaluations compared to the standard linear kernel approach, which relies on the use of simple inner product for the construction of Gramian matrix.

We focus on the bi-fidelity situation in this paper: We assume throughout well-defined low- and high-fidelity models exist, and that the low-fidelity model is inexpensive but also insufficiently accurate for the purpose of prediction. We assume that the HF model is an ideal computational simulation model that requires an onerous amount of time. Both LF and HF models are assumed to share a common set of tunable parameters (such as constitutive model constants or coefficients, and are often modeled as random variables for uncertainty quantification purposes). With a bi-fidelity template, the approaches in this paper can be extended to multi-fidelity \cite{zhu2014computational,razi2018MD}.

To demonstrate the effectiveness of our approach, we test the new procedures on several nontrivial examples:
\begin{enumerate}
  \item one- and two-species molecular dynamic (MD) simulation models: determination of MD model's statistical properties;
  \item polymer networks: investigation of the resultant strain of a graph of links in a polymer chain under a constant stress; 
\item $n$-body galaxy model: prediction of total energy and mean velocity and distance for the system; 
\item incompressible flow in two-dimensional channels with stenosis: estimation of velocity and wall shear stress profiles; and
  \item plasmonic nano-particle arrays: estimation of extinction and scattering efficiencies for Vogel spiral configurations with a different number of nano-particles.
\end{enumerate}
This manuscript is organized as follows: In the following two sections, the theoretical foundation of both multi-fidelity models and the proposed kernel selection approaches is provided. The sections are followed by the results and discussion of five canonical test problems. Finally, concluding remarks are presented as the final section of this paper, which is Section~\ref{sec:conclusion}. 
 
\section{Low-rank approximation using multiple levels of fidelity}\label{sec:MF}
Let a low-fidelity model $g_L$, and a high-fidelity model $g_H$, be given. These two models are functions of input parameters $p$. Consider a wave model with the wavespeed, $p$, as a parameter. For a given parameter value $p$, $g_L$ and $g_H$ output vectors:
\begin{align*}
  D &\subset \R^q, & g_L &: D \rightarrow \R^m, &
  g_H &: D \rightarrow \R^M.
\end{align*}
Note that usually $m \neq M$, and in general we assume $m \ll M$. (For example, a coarse $m$-degree-of-freedom model, versus a refined $M$-degree-of-freedom model.) We assume that $g_H(p)$ is a good predictor of some baseline ``truth'', but that $g_L(p)$ is a poor predictor. Qualitatively, $g_L$ is inexpensive but inaccurate, and $g_H$ is accurate but expensive.

The problem we consider is to construct an accurate emulator $\tilde{g}_H$ for $g_H$ whose accuracy is comparable to $g_H$ but whose cost is comparable to $g_L$, i.e., 
\begin{align*}
  \mathrm{Cost}(\tilde{g}_H) &\sim \mathrm{Cost}(g_L), & g_H(p) &\approx \tilde{g}_H(p), \hskip 10pt \forall\; p \in D.
\end{align*}
The focus of this paper is on such a construction of $\tilde{g}_H$. The end-user goal for this emulator may be for computation of statistics of $g_H$ (by modeling $p$ as a random variable), for determining sensitivity of $g_H$ to $p$, inference (identification of appropriate $p$ to match data). Any of the end-user goals require many queries of the $g_H$ for several values of $p$, which can be prohibitively expensive if $g_H$ is very expensive to query. This motivates the need to construct this emulator.

\subsection{Low-rank multi-fidelity model construction}\label{ssec:mf}
We describe the core idea behind low-rank MF methods, as outlined in \cite{narayan2014stochastic}. The construction of $\tilde{g}_H$ proceeds in two steps. In the first step, the parameter space $D$ is explored via the low-fidelity model to identify ``important'' points in parameter space. Since by assumption $g_L$ is an inexpensive emulator, this is practically feasible. The procedure for accomplishing this exploration is essentially numerical linear algebra.

First we densely sample $D$ with the ensemble $\{p_j\}_{j=1}^N$, where $N \gg 1$. These points may be selected in any reasonable way: as Monte Carlo samples, equispaced samples, low-discrepancy samples, etc. We next construct a Gramian matrix $G_L$ with entries:
\begin{align}
\label{eq:Eq1}
  (G_L)_{l,j} &= K\left(g_L(p_l), g_L(p_j)\right), &  l, j &= 1, \ldots, N, & K(u,v) &= \left\langle u, v \right\rangle,
\end{align}
where $\langle \cdot, \cdot\rangle$ is the standard Euclidean inner product on $\R^m$. This Gramian is positive semi-definite, but is not definite if $N < m$. In this case, the first $m$ steps of a pivoted Cholesky decomposition can be completed:
\begin{align}\label{eq:cholesky}
  P G_L P^T = L L^T,
\end{align}
where $L$ is a lower-triangular matrix and is zero for rows $m+1, \ldots, N$. The permutation matrix $P \in \R^{N \times N}$ determines an ``important'' point selection, in particular, the vector 
\begin{align*}
  z \coloneqq P \left(\begin{array}{c} 1 \\ 2 \\ \vdots \\ N \end{array}\right),
\end{align*}
is used to define an ordering $p_{z_1}, p_{z_2}, \ldots, p_{z_N}$ of the $N$ parameter points. For any $n \leq m$, the point set
\begin{align}
  \widehat{D} = \left\{ p_{z_1}, \ldots, p_{z_n} \right\} \subset D
\label{eq:Eq3}
\end{align}
are the ``important'' points. For later purposes we will also require an $n \times n$ matrix $\widehat{G}_L$, which is a $z$-sliced version of $G_L$:
\begin{align*}
  \left(\widehat{G}_L\right)_{i,j} &= \left(G_L\right)_{z_i, z_j}, & i, j &= 1, \ldots, n.
\end{align*}
The identification of $\widehat{D}$ is the main purpose of the first step. The second step will require evaluations of the HF model on $\widehat{D}$, so that $n$ is usually chosen based on the available computational budget. 

This approach for selecting ``important'' points via pivots is computationally efficient \cite{zhu2014computational,razi2018MD} and is equivalent to the pivoted samples selected by a column-pivoted QR decomposition of a matrix of model realizations~\cite{narayan2014stochastic}. Pivots from an LU factorization are also used in the area of computational fluid dynamics and heat transfer~\cite{anderson2017efficient} similarly. Alternative statistical strategies such as leverage-score sampling methods~\cite{perry2016augmented} and group matching methods \cite{lozano2011group,perry_allocation_2019} are also useful for the selection of optimal/important set of sampling points. For the purposes of simplicity in this paper, we focus exclusively on the Cholesky approach for identifying points.
The second step of the MF procedure first computes the HF ensemble,
\begin{align*}
  \left\{ g_H(p_{z_i}) \right\}_{i=1}^n,
\end{align*}
which is feasible if $n$ is small. The MF approximation is built as a linear expansion over this ensemble:
\begin{align}
  g_H(p) \approx \widetilde{g}_H(p) \coloneqq \sum_{l=1}^n g_H\left(p_{z_l}\right) c_l(p),
\label{eq:Eq4}
\end{align}
where $\{c_l(p)\}_{l=1}^n$ are coefficients in a vector $c(p)$ that is computed via a least-squares projection on the low-fidelity model:
\begin{eqnarray}\label{eq:Eq5}
  \widehat{G}_L c(p) 
    = \left( \begin{array}{c} K\left( g_L(p), g_L(p_{z_1}) \right) \\
                                     \vdots \\
    K\left(g_L(p), g_L(p_{z_n}) \right) \end{array}\right). \nonumber \\
\end{eqnarray}
Note that the right-hand side of the above equation is computable immediately simply by evaluating the $g_L(p)$, i.e., the low-fidelity model. \revision{The least squares projection coefficients can be obtained as
\begin{eqnarray}
c(p)=\widehat{G}_L^{-1} \left( \begin{array}{c} K\left( g_L(p), g_L(p_{z_1}) \right) \\
                                     
                                     \vdots \\
    K\left(g_L(p), g_L(p_{z_n}) \right) \end{array}\right). 
\label{eq:Eq005}
\end{eqnarray}
}
\rrevision{Note in particular that the emulator is interpolatory on the trained samples $\{p_i\}_{i=1}^n$, i.e., $g_H(p_i) = \widetilde{g}_H(p_i)$ for all $i = 1, \ldots, n$. This property can be directly deduced by noting in \eqref{eq:Eq5} that if $p = p_{z_i}$, then the right-hand side is exactly a column of $\widehat{G}_L$, and so all the $c_j(p_i)$ coefficients vanish, except for $c_j(p_j) = 1$.}
This procedure can be generalized to include more levels of fidelity~\cite{narayan2014stochastic,zhu2014computational}, which is advantageous when three or more models exist.

In summary, the first stage of the algorithm is ``offline'', which requires the computation of $N \gg 1$ realizations of the LF model. Some (inexpensive) linear algebra identifies a small set of $n \ll N$ parameter values where the HF model is evaluated in step 2. Then for each fixed $p \in D$, the low-rank MF surrogate in \eqref{eq:Eq4} is computed by first evaluating $g_L(p)$, and then by solving \eqref{eq:Eq5} for the vector $c(p)$.

While the MF surrogate $\widetilde{g}_H$ certainly achieves the cost requirement (evaluation requires only low-fidelity effort), understanding its accuracy is more subtle. The authors in \cite{narayan2014stochastic} provide abstract error certification in general Hilbert spaces; however the analysis is not useful in practice and the bound derived is not computable. In \cite{keshavarzzadeh_convergence_2019} a more practical analysis demonstrates that this surrogate is accurate when the LF and HF models correspond time integration solvers with coarse and fine timesteps, respectively. Finally, the authors in \cite{hampton2018practical} provide a statistical strategy for computationally certifying error committed by $\widetilde{g}_H$.

\subsection{Limitations and Drawbacks}\label{ssec:drawbacks}
Assuming $\widetilde{g}_H$ is accurate, the major limitations of this approach are described in Section \ref{sec:intro}. We give a more detailed discussion here. First recall that since $G_L$ has rank $m$ at most, then when $m \ll N$ only the first $n = m$ pivots selected in the decomposition \eqref{eq:cholesky} are informative. (The remaining $N-m$ pivots are essentially selected at random due to floating-point roundoff errors.) To see why this can become an issue, consider $g_L$ as a scalar output of a model (for example, the section lift coefficient in an aerodynamics model). Then the MF procedure only produces a single informative parameter value. This prevents the procedure from taking advantage of additional HF model runs when such runs are computationally feasible.

A second limitation is in the stability of the linear solver required in \eqref{eq:Eq5}. When $G_L$ is a low-rank matrix, then its singular values decay quickly, and the sliced matrix $\widehat{G}_L$ is ill-conditioned. This can numerically pollute the MF surrogate \eqref{eq:Eq4}.

A final limitation that we will discuss is more subtle: one requirement for success of this procedure is that $G_L$ is ``similar'' to $G_H$, where $G_H$ is a HF Gramian defined similarly as \eqref{eq:Eq1} but with $g_H$ instead of $g_L$. The choice of $K$ in \eqref{eq:Eq1} to take uncentered inner products may not reveal similarity of $G_L$ and $G_H$, but choosing a different function $K$ may improve this approximation.

\section{Kernel function Selection}\label{sec:kernel}
As discussed in Section~\ref{sec:MF}, kernel functions serve a key role in construction of the Gramian matrix, which is a building block of the low-rank multi-fidelity approximation approach. Equation \eqref{eq:Eq1} shows the appearance of the kernel $K$ in development of the low-rank MF model. Section \ref{sec:MF} shows the construction of the MF model with a linear kernel $K$ as shown in \eqref{eq:Eq1}. However, section \ref{ssec:drawbacks} lists some limitations of using this linear kernel. The main contribution of this paper is to mitigate some of these shortcomings by employing an alternative kernel function.

Note that there are several types of kernel functions that one could choose as an alternative to linear kernels. For simplicity, we consider several popular \textit{radial} kernels as alternatives. (A kernel $K$ is radial if $K(u,v) = F(\|u - v\|)$ for some function $F$.) We list these alternative kernels in Table \ref{tab:kernels}. Note that this list is not meant to be exhaustive, and only serves as a reasonably diverse set of kernels to demonstrate our approach. 

One a kernel is selected, say of index $i$, the entire procedure of Section \ref{sec:MF} can be completed, by replacing $K$ in \eqref{eq:Eq1} with $K_i$. (For example, the procedure remains unchanged if we consider $K_1$ since this is linear kernel in Table \ref{tab:kernels}.) There is one new detail in Table \ref{tab:kernels}, namely that most of the kernels are now functions of an additional set of hyperparameters $h$, where the number of hyperparameters depends on the family. The hyperparameters are frequently parameters that affect the general shape of the kernel, and it is well-known in both the computational mathematics and statistics community that a good choice for these parameters is essential to build both a stable and accurate model \cite{fasshauer_meshfree_2007,rasmussen_gaussian_2004}. Thus, we have complicated our procedure since we must now make additional choices for the family of kernel (quantified by the index $i$ in Table \ref{tab:kernels}), and we must choose values for the hyperparameters $h$. We first discuss the task of choosing values for the hyperparameters $h$.

\begin{table}
  \begin{center}
  \resizebox{\textwidth}{!}{
    \renewcommand{\tabcolsep}{0.4cm}
    \renewcommand{\arraystretch}{1.3}
    {\small
    \begin{tabular}{@{}lllp{0.8\textwidth}@{}}
      \toprule
      Index $i$ & $\dim h$ & Kernel type &  \\\midrule
      1 & 0 & Linear & $K_i = \left\langle u, v \right\rangle$ \\
      2 & 1 & Exponential & $K_i = \exp\left( - \frac{\|u - v\|}{h_1} \right)$\\
      3 & 1 & Squared exponential & $K_i = \exp\left( - \frac{\|u - v\|^2}{2 h_1} \right)$\\
      4 & 2 & Rational quadratic & $K_i = \left( \frac{1 + \|u - v\|^2}{2 h_1^2 h_2} \right)^{-h_2}$\\
      5 & 1 & Matern $3/2$ & $K_i = \left( 1 + \frac{\sqrt{3}}{h_1} \|u - v\|\right) \exp\left(-\frac{\sqrt{3}}{h_1} \| u - v\| \right)$\\
      6 & 1 & Matern $5/2$ & $K_i = \left( 1 + \frac{\sqrt{5}}{h_1} \|u - v\| + \frac{5}{3 h_1^2} \| u -v \|^2 \right) \exp\left(-\frac{\sqrt{5}}{h_1} \| u - v\| \right)$\\
      7 & 2 & Compact RBF & $K_i = \max \left[ 0, 1 - \left( \frac{1}{h_1} \| u - v\|\right)^{h_2} \exp \left( -\frac{\|u-v\|^2}{2 h_1^2} \right)\right]$\\
    \bottomrule
    \end{tabular}
  }
    \renewcommand{\arraystretch}{1}
    \renewcommand{\tabcolsep}{12pt}
  }
  \end{center}
  \caption{Kernel function library for $K(u,v;h)$. $u$ and $v$ are input vectors, while $h$ is a vector of hyperparameters. $\|\cdot\|$ denotes the Euclidean length of an input vector, and $\left\langle \cdot, \cdot \right\rangle$ denotes the Euclidean inner product. \revision{Kernel functions 2 through 7 are two types of radial kernel functions.}}\label{tab:kernels}
\end{table}

\subsection{Hyperparameter selection}\label{ssec:optimization}
In this section, let $i$ be a \emph{fixed} family index associated to one of the kernel families in Table \ref{tab:kernels}. Our task is to determine an appropriate value of the hyperparameter $h$. We formulate this selection as a simple optimization problem, where we impose that the following requirements should be met:
\begin{itemize}
  \item To maintain accuracy, the optimized kernel should behave like the linear kernel $K_1$ \revision{(since when the Gramian $G_L$ is well-conditioned, results in the literature \cite{narayan2014stochastic,keshavarzzadeh_convergence_2019} show that this procedure is effective)}. 
  \item To maintain stability, the Gramian $G_L$ resulting from a particular choice of kernel should reasonably be well-conditioned.
\end{itemize}
Our desire to maintain accuracy of the MF emulator is not always aligned with ensuring that a selected kernel behaves like the linear kernel. However, the current success of the linear kernel in the literature suggests that this may promote some accuracy\cite{razi2018quantized,razi2018MD,narayan2014stochastic}. In addition, the main purpose of the accuracy requirement is to ensure that the chosen hyperparameters do not force the kernel function to concentrate to a Dirac delta function: if the kernel focuses to the Dirac delta function, then the MF approximation will suffer poor accuracy. The choice of proximity to the linear kernel is thus just one choice of regularizer; we implement this as the Frobenius norm proximity to the linear kernel Gramian.

The need to maintain stability comes from the fact that solving \eqref{eq:Eq5} should be numerically stable. When the Gramian matrix is rank-deficient, the stability consideration regularizes the matrix so that its rank deficiency is mitigated. We implement this as the numerical (or stable) rank of the Gramian.

In summary, with $i \geq 2$ fixed, let $G_L^i(h)$ denote $N \times N$ the low-fidelity Gramian matrix associated with using $K = K_i(h)$ in \eqref{eq:Eq1}. I.e., $G_L^i(h)$ has entries
\begin{align*}
  \left(G_L^i(h)\right)_{j,k} = K_i( u_L(p_j), u_(p_k); h).
\end{align*}
We choose a value for the hyperparameters $h$ via the optimization problem
\begin{align}
  h^\ast = \argmin_h f_{\text{obj}}(h) \coloneqq \left\| G^1_L - G_L^i(h)\right\rVert_F + \lambda \frac{1}{\sqrt{\mathrm{srank}\left(G_L^i(h)\right)}} %\frac{\left\| G_L^i(h)\right\|_2}{\left\| G_L\right\|_F},
\label{eq:Eq6}
\end{align}
where $\lambda$ is tunable parameter, $\|\cdot\|_F$ is the Frobenius norm on matrices, and $\mathrm{srank}(\cdot)$ is the numerical/stable rank:
\begin{align*}
  \mathrm{srank}(A) = \frac{\|A\|^2_F}{\|A\|^2_2}.
\end{align*}
The objective function $f_{\textrm{obj}}$ in \eqref{eq:Eq6} thus seeks to simultaneously promote proximity to the linear Gramian (``accuracy") and maximization of the numerical rank (``stability").
The value of $\lambda$ can be obtained by empirical testing. 
\rrevision{In practice we judge that this procedure finds an effective kernel if $h^\ast \sim \lambda/\sqrt{n}$, corresponding to a case when the error in the adaptive kernel Gramian balances with the numerical rank. However, we expect optimal values of this parameter to be problem-dependent.}
\revision{For all problems tested in this paper, we found that $\lambda = 0.1$ provided good results, obtained through grid search and evaluating by computing errors on testing data. This value of $\lambda$ was computed using a logarithmic grid search over the first three test problems in this manuscript. Using this approach, $\lambda$ is a regularization parameter that promotes the numerical stability of the solution, and the grid search we have employed is one of the standard approaches in machine learning~\cite{bergstra2012random,klein2017fast} for tuning such a regularization parameter.}

As often the number of hyperparameters is small (cf. Table \ref{tab:kernels}) and the computational cost of high-fidelity model evaluation is significantly higher than  that of the low-fidelity model, the process of hyperparameter optimization is inexpensive in contrast. Note that the optimization \eqref{eq:Eq6} operates \textit{only} on the low-fidelity model, and the objective function is simple and efficient to evaluate if all the low-fidelity data is stored.

In order to deal with this problem, we propose two approaches that generalizes the kernel function selection framework. For both of these proposed approaches we are using a heuristic optimization algorithm called particle swarm algorithm~\cite{kennedy2011particle,clerc2010particle} to minimize $f_{\text{obj}}$ in Eq.~\ref{eq:Eq6}. Next, this algorithm is briefly discussed.

\subsection{Particle swarm optimization}\label{sec:PSO}
The optimization problem \eqref{eq:Eq6} has $\dim h$ design variables, and for our choices of kernels this number is relatively small (see Table \ref{tab:kernels}). Nevertheless, gradient-based optimization algorithms can perform poorly when the initial guess is far away from a basin of attraction. To mitigate this issue, we employ particle swarm optimization \cite{kennedy2011particle,clerc2010particle} to compute a reasonable initial guess, and subsequently perform a gradient-based interior point algorithm. 

Particles swarm optimization (PSO) is a velocity-based evolutionary optimization approach, and was initially developed for modeling social behaviors in bird flocks using topological neighborhoods~\cite{deepa2011model}. PSO is initialized using an ensemble of random points in the design/search space, with each ensemble member called a particle. For each particle there is an assigned position (design variable value) and velocity. The initial selection of particle positions and velocities is often random from a uniform distribution. The remaining portion of the algorithm updates particle positions and velocities using a notion of \textit{fitness} of a particle, i.e., an estimate of its optimality relative to the ensemble.

Let $\{h_i\}_{i = 1}^S$ and $\{v_i\}_{i=1}^S$ be the position and velocity, respectively, of the $i$th particle. For $n \in \N$, let $h_i^n$ denote the position of the $i$th particle at PSO iteration $n$. At $n = 1$, all positions and velocities are randomly initialized.
PSO updates proceed by:
\begin{eqnarray}\label{eq:Eq7}
  v^{n+1}_i&=&v^n_i+k_1\rho \left(B^n_i-h^n_i\right)+k_2\gamma \left(B^n_g-h^n_i\right), \\ \nonumber 
  h^{n+1}_i &=& h^n_i+v^{n+1}_i,
\end{eqnarray}
where $B^n_i$ and $B^n_g$ denote the best solution for the $i$th particle in iterations $1,...,n$, and group's best solution for iterations $1, \ldots, n$, respectively. I.e., 
\begin{align*}
  B_i^n &\coloneqq \min_{j \in [n]} f_{\text{obj}}(h_i^j), & B_g^n &\coloneqq \min_{(i,j) \in [S] \times [n]} f_{\textrm{obj}}(h_i^j), & [n] &\coloneqq \{1, 2, \ldots n\}.
\end{align*}
The parameters $\rho$ and $\gamma$ in \eqref{eq:Eq7} are uniform random variables on $[0,1]$ (independent at each time and for each particle) and $k_1$ and $k_2$ are tunable parameters. We apply standard velocity clamping procedures to ensure velocities remain bounded and particles do not leave the design space: If some particle's velocity exceeds a maximum allowed velocity limit, it is set to the maximum allowed value of velocity~\cite{shahzad2009opposition}.
PSO iterations continue until a stopping criterion is triggered. In our case, we limit the maximum number of iterations and lack of improvement in the group's best solution. 
There are many variations of PSO; for the purpose of this study we employ the MATLAB toolbox ``particleswarm solver''. This algorithm has been shown to be faster than standard evolutionary algorithms in approaching the vicinity of the minima~\cite{angeline1998using,angeline1998evolutionary}. 

However, PSO is slow when it comes to converging to a local minimum once a basin of attraction is reached. Therefore, we use loose criteria to terminate the PSO iterations, and subsequently use the best particle as the initial guess into a gradient-based interior point algorithm (with finite difference approximated gradient); the output of the gradient-based algorithm is our computed value of $h^\ast$ in \eqref{eq:Eq6}.

\subsection{Additive kernel approach}\label{sec:add}
We have discussed in Sections \ref{ssec:optimization} and \ref{sec:PSO} how we select the hyperparameter value $h$ for each kernel index $i$ in Table \ref{tab:kernels}. This section and the following section propose two strategies for choosing an overall kernel $K$ for use in the MF procedure of Section \ref{ssec:mf}.

Our first approach, the \emph{Additive Kernel Approach}, choose an overall kernel $K$ as an additive mixture of all kernels in the library:
\begin{align}
  K(u,v)&=\sum_{i=1}^L w_i K_i(u,v; h_i^\ast), & (w_1, \ldots w_L) &\in \{ x \in \R^L\; \big|\; \sum_{\ell=1}^L x_\ell = 1, \;\; x_\ell \in [0,1] \;\; \forall i\}
\label{eq:Eq8}
\end{align}
where $L$, $h^\ast_i$ and $w_i$ are the number of kernel functions in the library ($L=7$ in this paper), the $i$th optimized hyperparameter value $h^\ast$ output from \eqref{eq:Eq6}, and its corresponding weight. 
The optimizations associated with computing the quantities in \eqref{eq:Eq8} are convex, so several types of algorithms are appropriate. We again employ PSO for an initial guess, followed by an interior-point algorithm for accelerated convergence to a minimum.

The full MF procedure is exactly as described in Section \ref{ssec:mf}, with $K$ in \eqref{eq:Eq1} the kernel computed in \eqref{eq:Eq8}. We emphasize again that the choice of kernel \emph{only} depends on low-fidelity data. Only \emph{after} a kernel function is chosen do we query the HF model through the procedure in Section \ref{ssec:mf}.

\subsection{Adaptive approach for the selection of an optimal kernel function}\label{sec:adaptive}
Our second procedure for selecting a kernel is the \textit{Adaptive Kernel Approach}. In this procedure, we select \emph{only} one kernel $K_i$ as the kernel of choice in the MF procedure. However, the index of this choice $i$ depends on the HF simulation budget $n$. In brief, with $n$ fixed, we select $i$ as the kernel index that results in the best approximation of a size-$n$ least squares approximation of the LF model.

In detail, we first fix $i$ and construct $\widetilde{g}_{L,i}$ as an emulator for the LF model using kernel $i$:
\begin{align}
  \widetilde{g}_L(p)=\sum_{l=1}^n g_L\left(p_{z_{l,i}}\right) c_{l,i}(p),
\label{eq:Eq9}
\end{align}
where $\{z_{l,i}\}_{l=1}^n$ are the pivots selected from the Cholesky decomposition of $G_L$ in \eqref{eq:Eq1} with $K = K(\cdot, \cdot; h^\ast_i)$. The coefficients $\{c_{l,i}\}_{l=1}^n$ are entries in the least-squares coefficient vector $c_i(p)$ from a linear system similar to \eqref{eq:Eq5}:
\begin{eqnarray}
  \widehat{G}_L c_i(p) 
    = \left( \begin{array}{c} K\left( g_L(p), g_L(p_{z_{1,i}}) \right) \\
                                     \vdots \\
    K\left(g_L(p), g_L(p_{z_{n,i}}) \right) \end{array}\right), \nonumber
%\label{eq:Eq10}
\end{eqnarray}
where again $K = K(\cdot,\cdot; h^\ast_i)$. The function $\widetilde{g}_L$ is an emulator for $g_L$ built via a least squares approximation. One error metric for this least squares emulator is
\begin{align*}
  \epsilon_i \coloneqq \textrm{med}_{j \in [N] \backslash \{z_{1,i}, \ldots, z_{n,i}\}} \left\| g_L(p_j) - \widetilde{g}_L(p_j) \right\|,
\end{align*}
where $\textrm{med}$ denotes the median. Note that the error $\left\| g_L(p_j) - \widetilde{g}_L(p_j) \right\| = 0$ when $j \in \{p_{z_{1,i}}, \ldots, p_{z_{n,i}}\}$ represents the indices of the training samples. \revision{The median has been chosen as the summary statistic for the distribution of the residuals for multi-fidelity emulators since  the errors do not necessarily follow a normal (or symmetric) distribution and it is most often skewed. 
Moreover, the median is an effective statistic for summarizing distributions~\cite{dittmann2008biases} and thus provides a more robust measure of error for our purpose.}
\rrevision{We emphasize that the data used to train the emulator is the set of parameters $\{p_{z_i}\}_{i=1}^n$. In this testing phase we compute errors over the set of realizations corresponding to the parameters $\{z_i\}_{i=1}^N \backslash \{p_{z_i}\}_{i=1}^n$ that the emulator was \textit{not} trained with.}

To choose the ``best'' kernel index $i^\ast$, we choose
\begin{align*}
  i^\ast = \argmin_{i = 1, \ldots, n} \epsilon_i.
\end{align*}
Finally, we choose $K$ in \eqref{eq:Eq1} as $K = K_{i^\ast}(h^\ast)$. Like the Additive kernel approach, this approach analyzes \emph{only} the LF model, and hence is computationally feasible. However, this approaches differs from the additive approach in that the optimal $i^\ast$ is recomputed when $n$ is changed.

\section{Results and Discussion}\label{sec:results}
\newcommand{\angstrom}{\textup{\AA}}
In order to demonstrate the two proposed kernel function selection methods for constructing an accurate surrogate, we consider their implementation on five nontrivial technical problems from different areas of engineering and science. These problems have been selected to provide intuition about the feasibility and benefits of the application of both approaches to real-world problems with low-dimensional output parameter space. The success of these approaches can be directly measured by the accuracy of their MF surrogate model. Our chosen error metric is the median relative error:
  \begin{equation}
    \mathrm{error} = \underset{p \in \{p_1, \ldots, p_N\} \backslash \{p_{z_1}, \ldots, p_{z_n}\}}{\mathrm{median}} \frac{\left\| g_H(p) - \widetilde{g}_{H}(p) \right\|_{\R^M}}{\|g_H(p)\|_{\R^M}},
\label{eq:F1}
  \end{equation}
where $\{z_1, \ldots, z_n\}$ are the indices where the HF model was evaluated to construct $\widetilde{g}_H$.
Thus, our metric is evaluated over data from full size-$N$ discretized parameter set, and uses the high-fidelity model $g_H$ as the oracle ``truth''. We also compare results obtained with the use of the linear kernel function to emphasize the enhancement provided by the Additive kernel and Adaptive kernel approaches. We also measure cost of all procedures in units of the cost of a single HF realization. We include the cost of optimization in the kernel selection process as well. Thus, our cost metric is the ``effective number of high-fidelity samples", defined as
\DeclarePairedDelimiter\ceil{\lceil}{\rceil}
\DeclarePairedDelimiter\floor{\lfloor}{\rfloor}
  \begin{eqnarray}
    \text{Effective Number of High-fidelity Samples}=\text{Number of HF Samples} \nonumber \\
    \mbox{\ } +\ceil*{\frac{\text{Cost(Kernel optimization)}}{\text{Cost(One HF Simulation)}}},
\label{eq:F2}
  \end{eqnarray}
where $\ceil*{}$ is the ceiling function. Many of our plots will measure accuracy versus cost, using the accuracy and cost metrics identified above. 

\subsection{Test Problem 1: Molecular systems}

The classical computational algorithm for performing molecular dynamics (MD) simulations has been developed based upon the Lagrangian methodology of tracking particle dynamics via Newton's equations of motion. Numerical integration of these equations provides an accurate evaluation of the time evolution of a molecular system, and hence the system's quantities of interest. While there are many integration schemes available in the literature, but we settle on the velocity Verlet algorithm~\cite{lee2016computational} due to its popularity, reasonable accuracy, and simplicity of implementation. The Lennard-Jones (LJ) potential was used for simulating the interatomic interactions in the MD system. The potential function also has parameters $\epsilon$, the potential well depth, and $\sigma$, the length scale for the pairwise interatomic interaction.

In the process of MD simulation using any interatomic interaction model, when the size of the integration time-step is large, the major challenge is the stability of the MD integration scheme. If the time-step is too large, then iterations may predict near-collocation of molecules/atoms which will lead to large (unphysical) repulsive forces and large displacements on the next time step, and this instability snowballs to future time steps. To prevent such divergence of the integration scheme, we clamp the magnitude of repulsive interactions for closely-approaching atoms. For our LJ potential, this capping can be implemented straightforwardly by modifying the potential at short distance.

The MF situation we consider here are models whose fidelity is defined by the choice of integration time step. The parameters $p \in \R^2$ we consider are the temperature and the density of the molecular system. For every parameter value, quantities of interest are computed using the low-fidelity model (here, the MD simulation with a large time-step $\Delta t$) and concatenated into an output vector. The quantities of interest we consider here are the scalar averaged total energy $E$ and diffusion coefficient $D$:
\begin{equation}
  g(p_j)=\left[\frac{E(p_j)}{\sqrt{\overline{E}}},\frac{D(p_j)}{\sqrt{\overline{D}}}\right]^T.
\label{eq:P1_1}
\end{equation}
We have rescaled the parameters by their averaged values over the ensemble:
\begin{eqnarray}
 \overline{E} &=& \frac{1}{N} \sum_{j=1}^N \| E(p_j) \|_2^2,  \nonumber\\
 \overline{D} &=& \frac{1}{N} \sum_{j=1}^N \| D(p_j) \|_2^2.
\label{eq:P1_2} 
\end{eqnarray}
Note that in MD simulations there are other important quantities of interest, such as the radial distribution function (RDF) and mean squared displacement (MSD).
For this test problem, we consider a uniform grid of temperature and density defined as $T\times \rho: [500, 1000]\text{K}\times[36.27,701.29]\frac{\text{kg}}{\text{m}^3}$ ($\rho^{\ast}: [0.05,0.95]$) with $N=114$ sample points, where  $\rho$ and $\rho^{\ast}$ are density and dimensionless density equal to $\frac{Nm}{V}$ and $\frac{N\sigma^3}{V}$, respectively. We also assume $\sigma$, $\varepsilon$, $m$ or molecular mass to be $3\angstrom$, $1\frac{\text{kcal}}{\text{mol}}$ and 12.01$\frac{\text{g}}{\text{mol}}$, respectively. Here, the boundary conditions  in all sides of cubic simulation box (with the width of 27.05$\angstrom$) are considered to be periodic. Here the mass of each particle ($m$) and simulation box size ($L=V^{\frac{1}{3}}$) is set to 12.01 g/mol and 27.05$\angstrom$, respectively. 

\begin{figure}[!tb]
  \centering
  \begin{subfigure}[t]{0.4\textwidth}
    \centering
    \includegraphics[width=\textwidth]{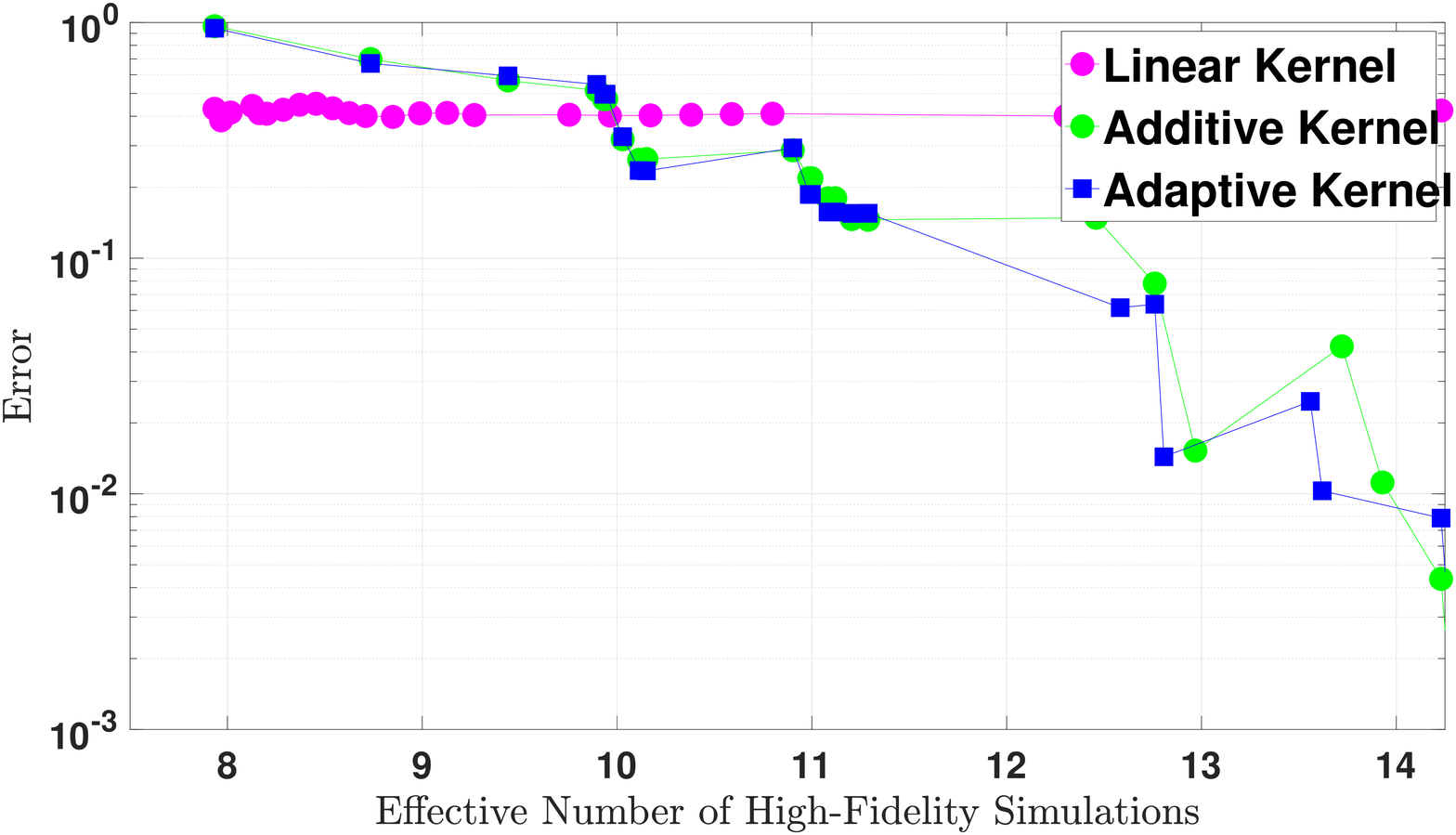}
    \caption{Total Energy}
  \end{subfigure}
  %\subfloat[Total Energy]{\label{Afig1}\epsfig{figure=figure/MD1_E.eps,width=0.48\textwidth}} \quad
  ~\begin{subfigure}[t]{0.4\textwidth}
    \centering
    \includegraphics[width=\textwidth]{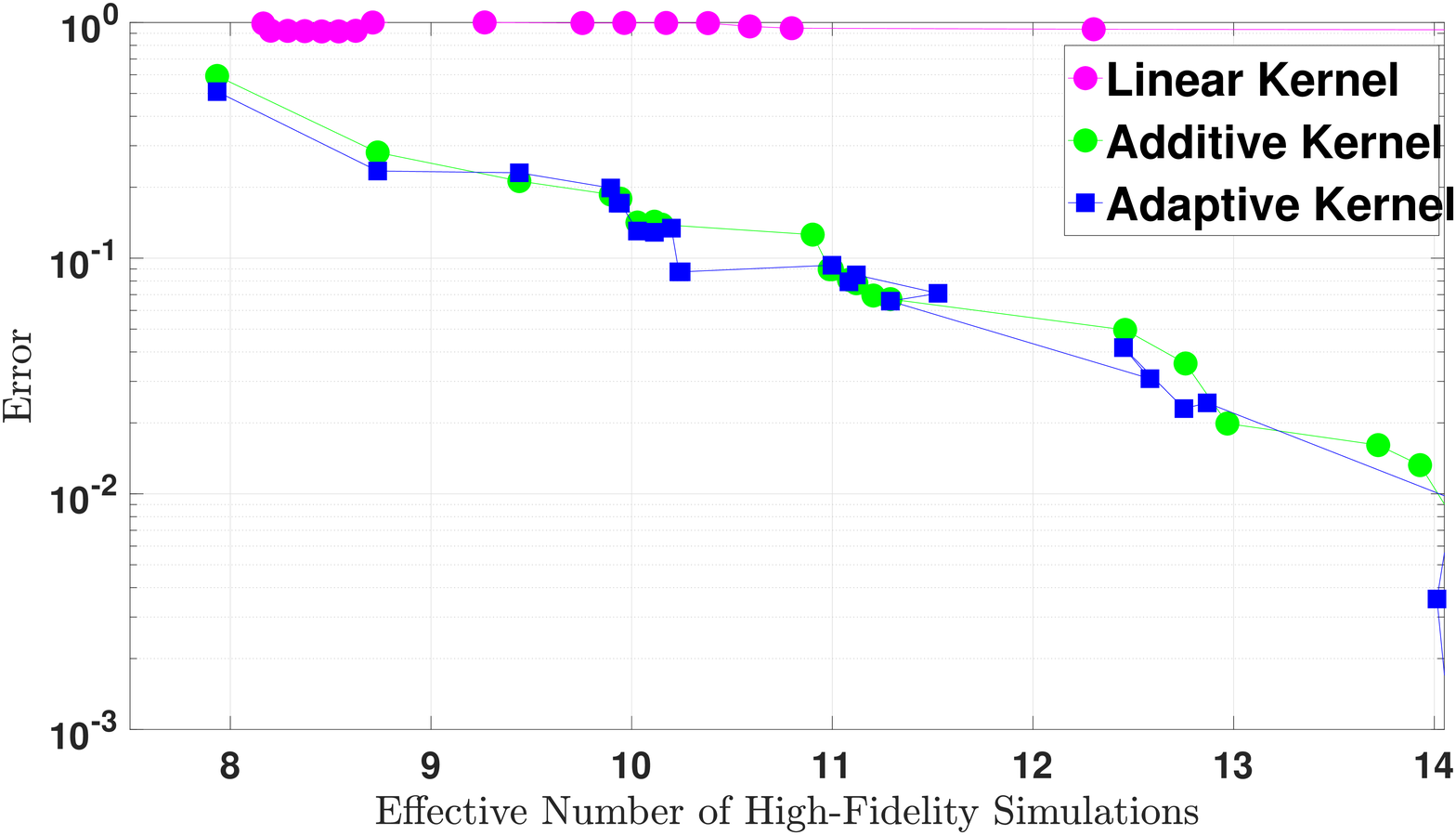}
    \caption{RDF}
  \end{subfigure}\\
  %\subfloat[RDF]{\label{Bfig1}\epsfig{figure=figure/MD1_RDF.eps,width=0.48\textwidth}} \newline\newline\newline\vskip -0.5cm
  ~\begin{subfigure}[t]{0.4\textwidth}
    \centering
    \includegraphics[width=\textwidth]{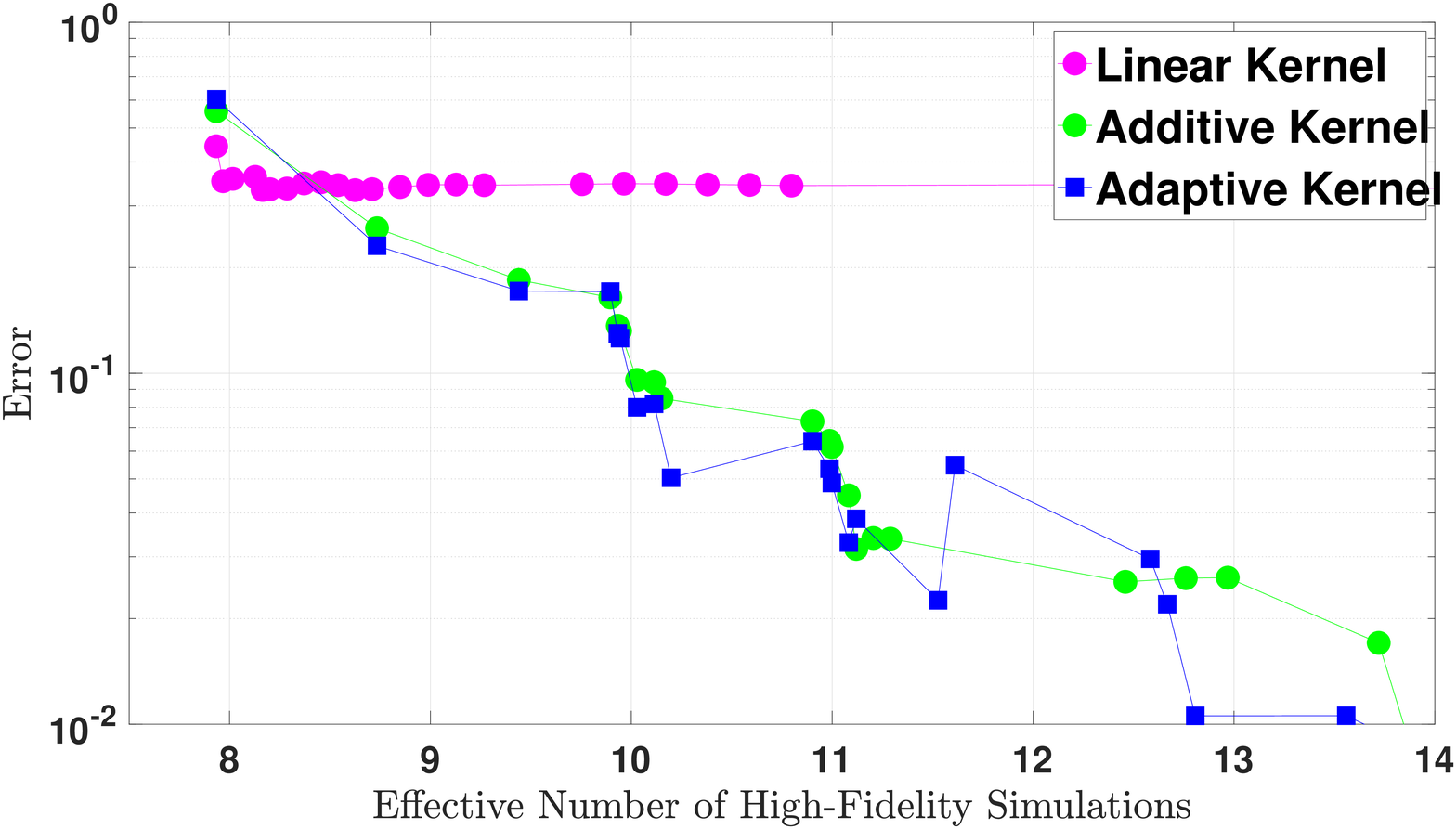}
    \caption{Self Diffusion Coefficient}
  \end{subfigure}
%\subfloat[Self Diffusion Coefficient]{\label{Cfig1}\epsfig{figure=figure/MD1_D.eps,width=0.48\textwidth}} \quad
  ~\begin{subfigure}[t]{0.4\textwidth}
    \centering
    \includegraphics[width=\textwidth]{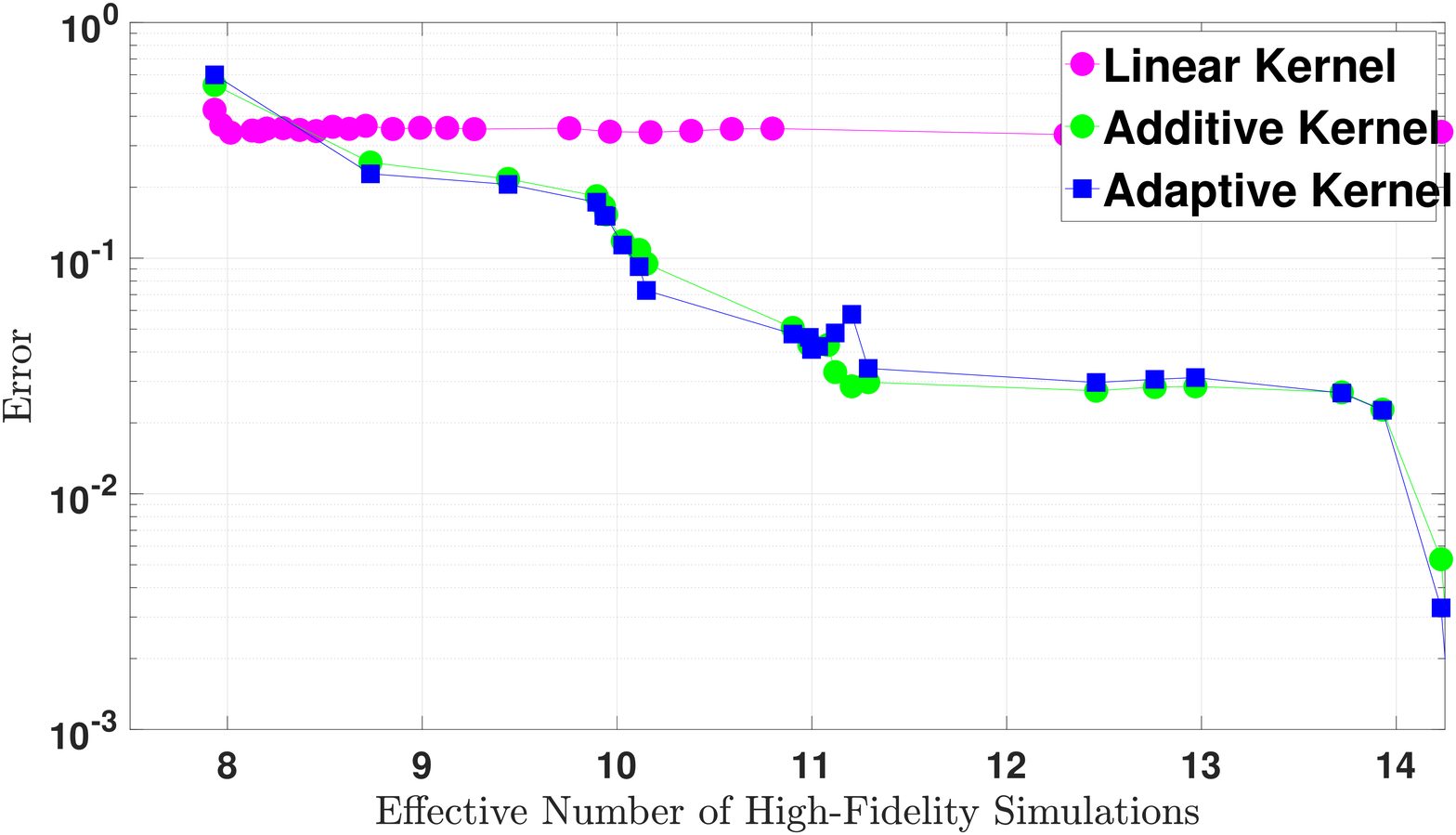}
    \caption{MSD}
  \end{subfigure}
%\subfloat[MSD]{\label{Dfig1}\epsfig{figure=figure/MD1_MSD.eps,width=0.48\textwidth}} \quad
\caption{Median error of the multi-fidelity model constructed based upon the results from the low rank low-fidelity MD model for one-species LJ data with $\Delta t=20$fs in the prediction of the properties for the high-fidelity model with $\Delta t=1$fs; Test Problem 1.}\label{fig1}
\end{figure}

The LF model is defined by an MD simulation with time step $\Delta t = 20$fs, and the HF model is the same MD simulation with $\Delta t = 1$fs. In this situation, the dimension of the vector $g_L(p)$ is 2, and hence the naive MF procedure in Section \ref{ssec:mf} with the linear kernel can only take $n = 2$ meaningful HF samples. Figure \ref{fig1} demonstrates that the Additive/Adaptive kernel approaches circumvent this limitation. We see that the linear kernel error stagnates since after $n = 2$ the ``important'' parameter values are chosen essentially randomly, and furthermore the least squares problem \eqref{eq:Eq5} becomes ill-conditioned. However, both the Additive and Adaptive kernel approaches produce more reasonable errors, with improvements gained as one is able to invest more HF effort.

As shown in Fig~\ref{fig1} (a) and (c), both proposed approaches work well although the error is not monotonically decreasing with respect to increased effort. \revision{This lack of monotonic error is due to the fact that we use kernel evaluations of the low-fidelity model on the right-hand side of \eqref{eq:Eq5} in the place of the (unavailable) high-fidelity kernel evaluation.} In panes (b) and (d) we show the error of the RDF and MSD quantities of interest for this procedure, which were \emph{not} used in formation of the MF surrogate. Thus, this procedure does not overfit energy and diffusion at the cost of other quantities.

Another reason for the non-monotonic reduction of the error with respect to the increase in effective high fidelity simulation is variable computational costs of molecular simulation, which are directly proportional to molecular density or the number of particles in the simulation box. Hence, addition of one high-fidelity sample does not linearly increase the cumulative computational cost.

\revision{
Here, both proposed approaches outperform the naive multi-fidelity method that applies the linear kernel function. The rank deficiency of the constructed Gramian matrix for the linear kernel approach is the main reason for such poor performance. This ill-conditioned Gramian matrix itself results from the low dimension of the model output space for low- and high-fidelity models.

Additionally, the results illustrated in Fig~\ref{fig1} indicate about two orders of magnitude in error reduction when the effective number of high-fidelity simulation increase from 8 to 14. This is supporting evidence for the success of the proposed multi-fidelity kernel optimization approach in handling a difficult predictive task for a complex system with underlying stochastic behavior.
}

Next, we consider a more difficult MD simulation with a two-component glass-forming system. In this case, the three output parameters of both LF and HF models are the diffusion coefficients $D_A$ and $D_B$ for species A and B, along with the total energy $E$.

\begin{figure}[!tb]
  \centering
  \begin{subfigure}[t]{0.4\textwidth}
    \centering
    \includegraphics[width=\textwidth]{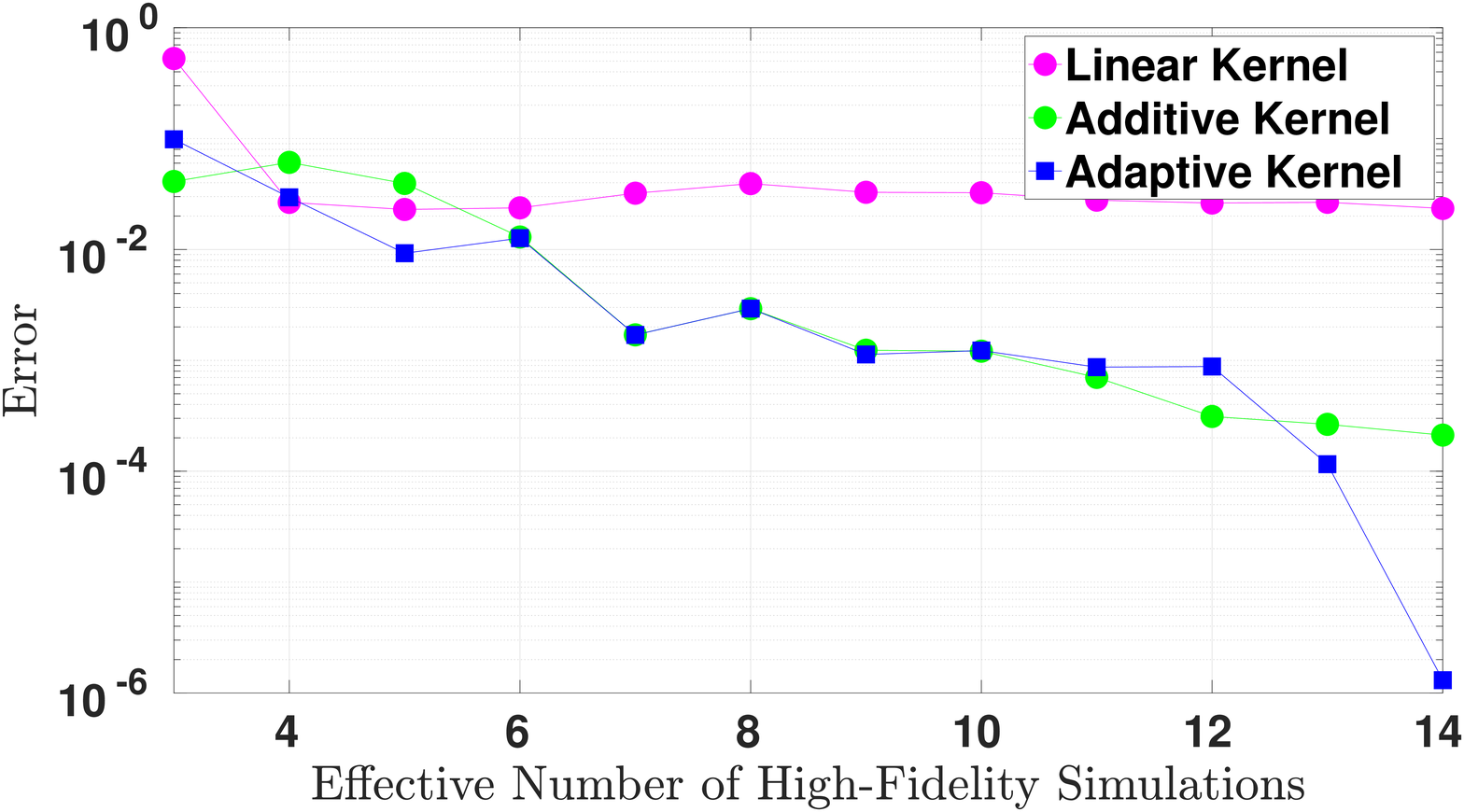}
    \caption{Total Energy}
  \end{subfigure}
  %\subfloat[Total Energy]{\label{Afig2}\epsfig{figure=figure/MD2_E.eps,width=0.48\textwidth}} \quad
  ~\begin{subfigure}[t]{0.4\textwidth}
    \centering
    \includegraphics[width=\textwidth]{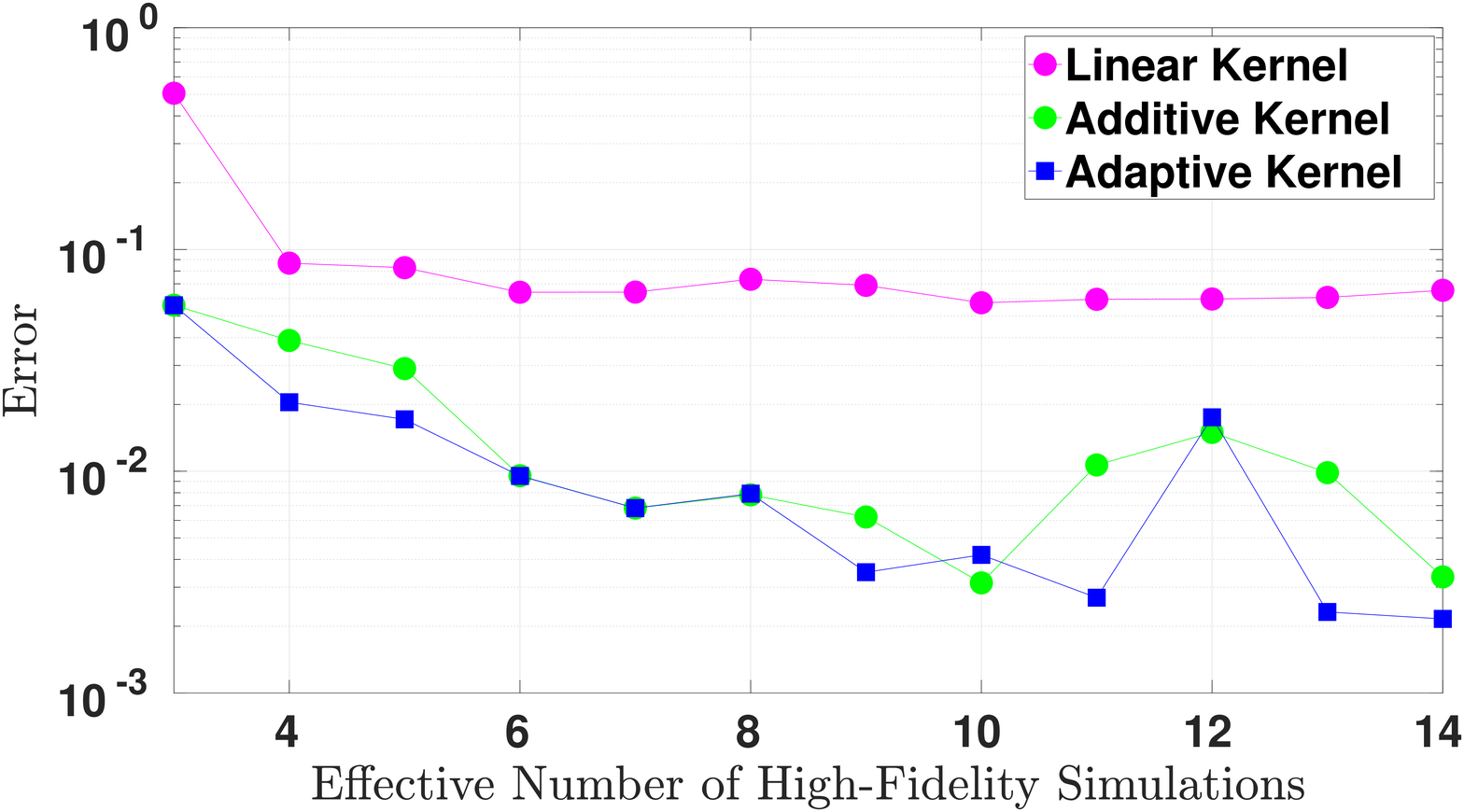}
    \caption{RDF A-A interaction}
  \end{subfigure}\\
  %\subfloat[RDF A-A interaction]{\label{Bfig2}\epsfig{figure=figure/MD2_UAA.eps,width=0.48\textwidth}} \newline\newline\newline\vskip -0.5cm
  \begin{subfigure}[t]{0.4\textwidth}
    \centering
    \includegraphics[width=\textwidth]{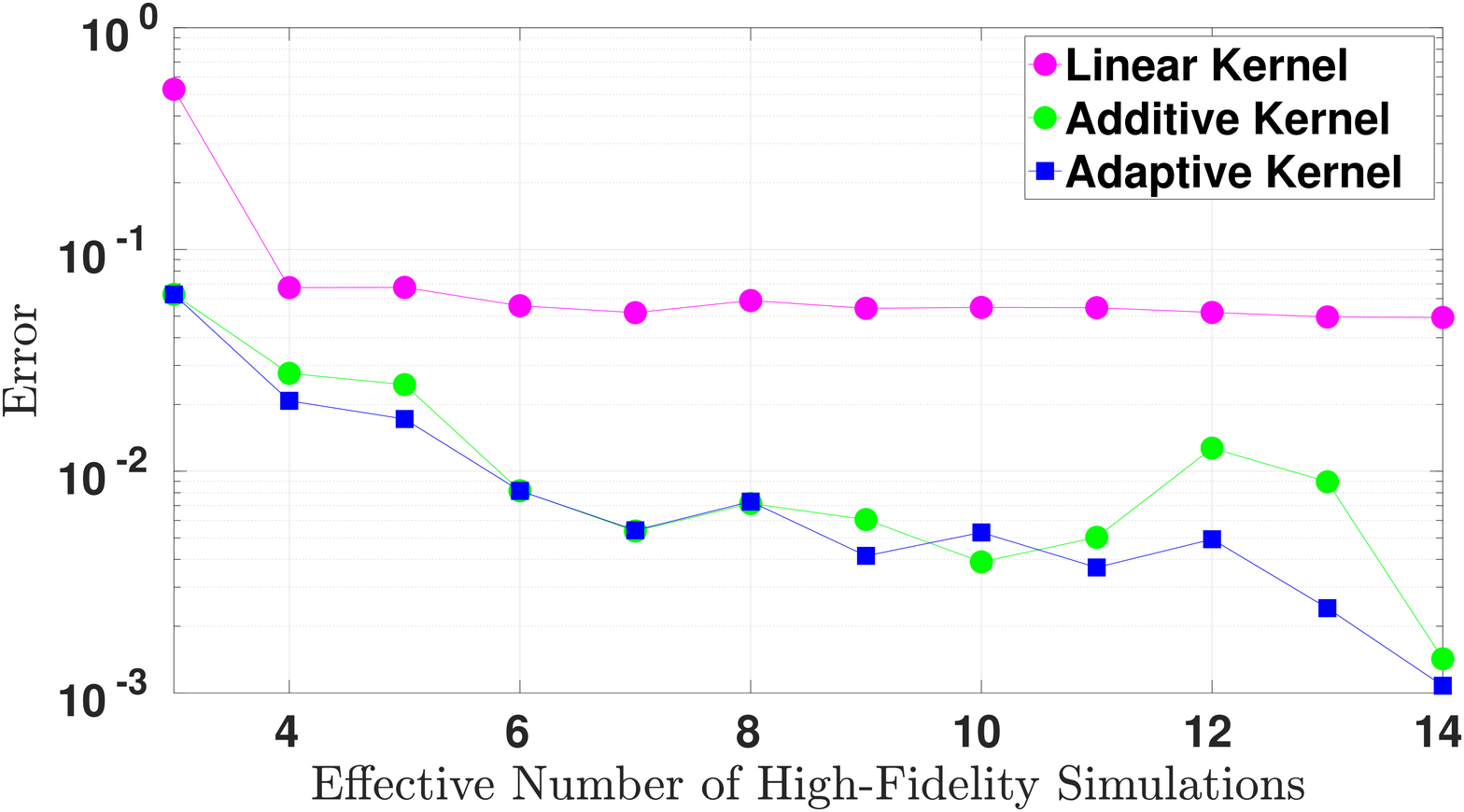}
    \caption{RDF A-B interaction}
  \end{subfigure}
  %\subfloat[RDF A-B interaction]{\label{Cfig2}\epsfig{figure=figure/MD2_UAB.eps,width=0.48\textwidth}} \quad
  ~\begin{subfigure}[t]{0.4\textwidth}
    \centering
    \includegraphics[width=\textwidth]{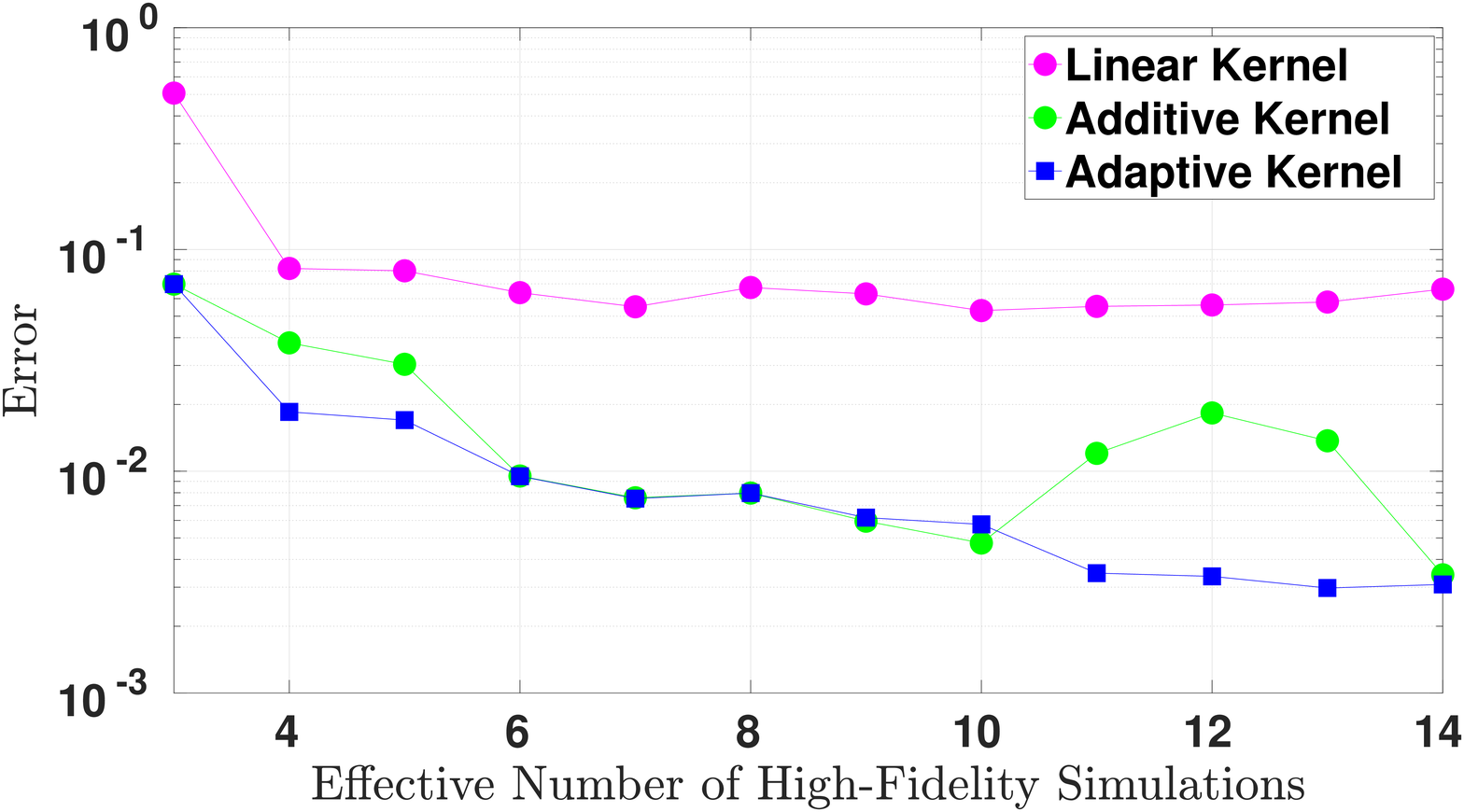}
    \caption{RDF B-B interaction}
  \end{subfigure}\\
  %\subfloat[RDF B-B interaction]{\label{Dfig2}\epsfig{figure=figure/MD2_UBB.eps,width=0.48\textwidth}} \newline\newline\newline\vskip -0.5cm
  \begin{subfigure}[t]{0.4\textwidth}
    \centering
    \includegraphics[width=\textwidth]{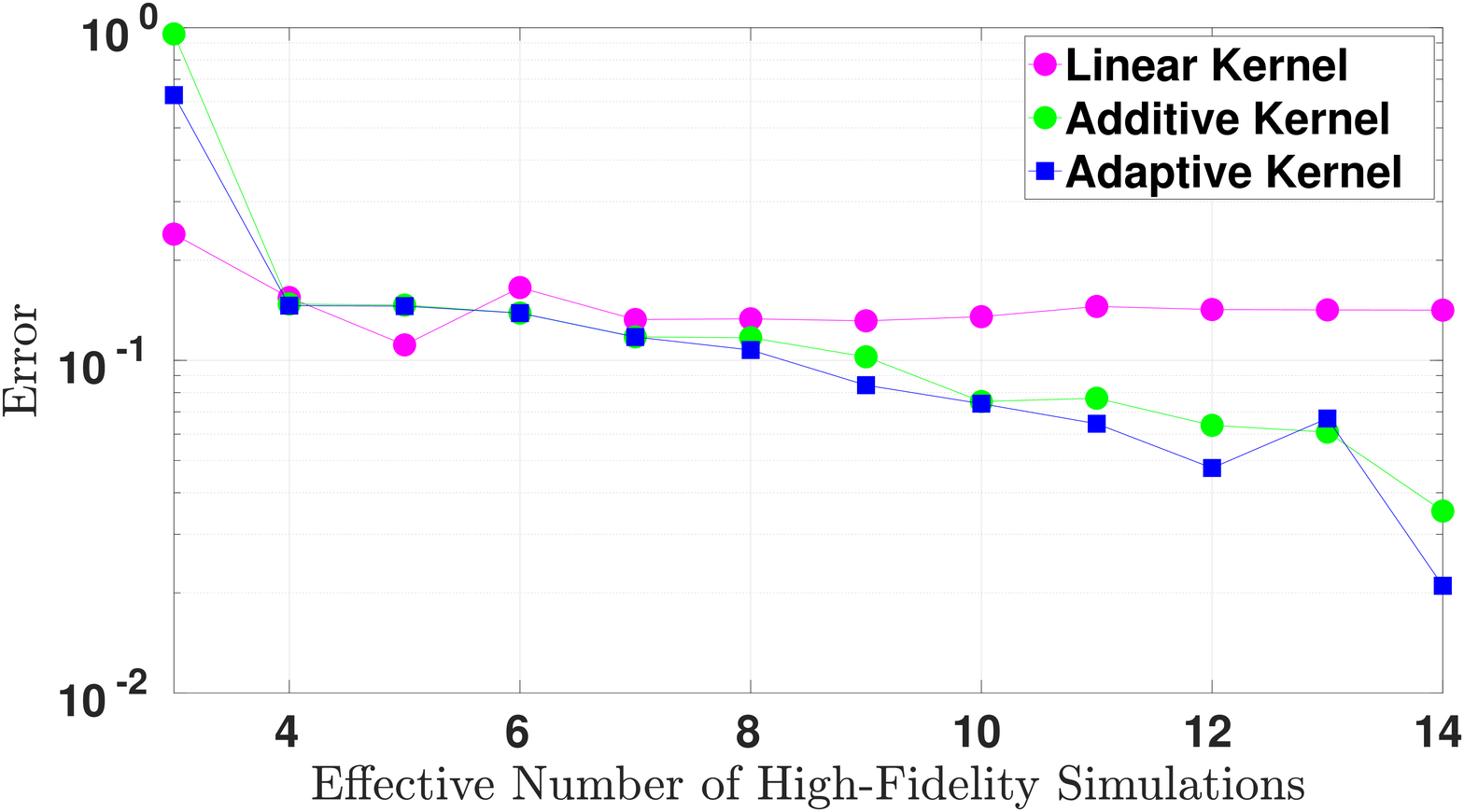}
    \caption{Self Diffusion Coefficient for species A}
  \end{subfigure}
%\subfloat[Self Diffusion Coefficient for species A]{\label{Efig2}\epsfig{figure=figure/MD2_DA.eps,width=0.48\textwidth}} \quad
  ~\begin{subfigure}[t]{0.4\textwidth}
    \centering
    \includegraphics[width=\textwidth]{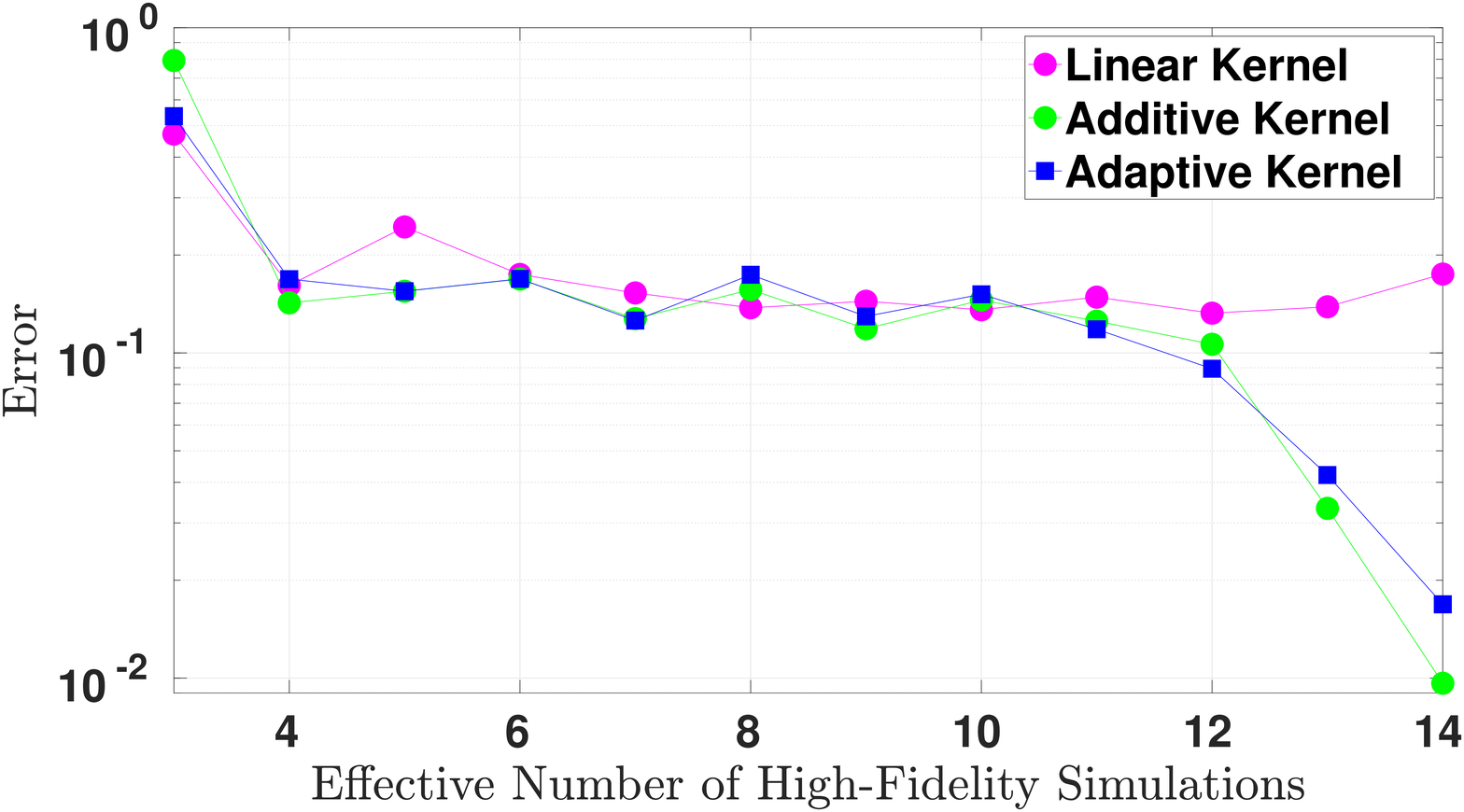}
    \caption{Self Diffusion Coefficient for species B}
  \end{subfigure}
%\subfloat[Self Diffusion Coefficient for species B]{\label{Ffig2}\epsfig{figure=figure/MD2_DB.eps,width=0.48\textwidth}} \quad
\caption{Median error of the multi-fidelity model constructed based upon the results from the low rank low-fidelity MD model for two-species LJ data with $\Delta t=20$fs in the prediction of the properties for the high-fidelity model with $\Delta t=1$fs; Test Problem 1.}
\label{fig2}
\end{figure}

We perform output normalization similar to \eqref{eq:P1_2}. Here, the length scale for pairwise interatomic interactions between molecules are set as $\sigma_A=3.5\angstrom$ for interactions between molecules of species A, $\sigma_B=3.888\angstrom$ for B-B interactions, and $\sigma_{AB}=\sqrt{\sigma_A\sigma_B}$ for A-B interactions interaction between both types of molecules. For this experiment, we have $N=24$ points in the density-temperature (parameter) space. The number of molecules in the MD simulation for each species is 512. Once again, the boundary conditions are assumed to be periodic for all sides of the simulation box. The temperature and simulation box length for these sample points vary between 180K and 290K, and 37.62$\angstrom$ and 38.35$\angstrom$, respectively. Once again, the results depicted in Fig.~\ref{fig2} shows high accuracy of both adaptive and additive models for prediction of the three radial distribution functions involving interactions of molecules of type A and B ($R_{AA}$, $R_{AB}$ and $R_{BB}$) as well as the corresponding diffusion coefficients and system total energy. \revision{However these results indicate that, if performing 12 high-fidelity simulation is computationally affordable, the multi-fidelity model constructed based on the proposed adaptive kernel function selection approach is slightly more accurate compared to the one constructed using the additive approach and significantly more accurate when compared to the multi-fidelity solution with linear kernel function.}

\revision{
Furthermore, as the number of particles is fixed for each sample, a more monotonic trend of convergence can be observed for this case (see Fig.~\ref{fig2}). Due to the smaller number of training low-fidelity data points, the accuracy of the approach degrades. On the other hand, while radial distribution functions are not involved in the construction of the multi-fidelity models, the accuracy of the emulator for their estimation is quite good: around 1 percent error with only 1/3 of the high-fidelity samples.
}

\subsection{Test Problem 2: The $n$-body galaxy model} 
Galaxy systems often can be modeled as an $n$-body problem, in which the gravitational dynamics governs the interaction between objects. The interaction force between the objects can be driven from Newton's law of gravity as~\cite{trenti2008gravitational}
\begin{equation}
F_{i}=-\sum_{\substack{j=1\\j \neq i}}^{n}{\frac{Gm_im_j\left(r_i-r_j\right)}{|r_i-r_j|^3}} ,
\label{eq:P2_1}
\end{equation} 
where $G$, $m_i$ and $r_i$ are the gravitational constant, $i$th object mass, and distance of $i$th object from the origin, respectively. Considering a cutoff radius can simplify the calculation, as the force between each object and any other object outside this radius is automatically set to zero. For the purpose of this study, the GalaxSee $n$-body model~\cite{jacobs2014blue} is used as the simulation module. 

\begin{figure}[!tb]
  \centering
  \begin{subfigure}[t]{0.4\textwidth}
    \centering
    \includegraphics[width=\textwidth]{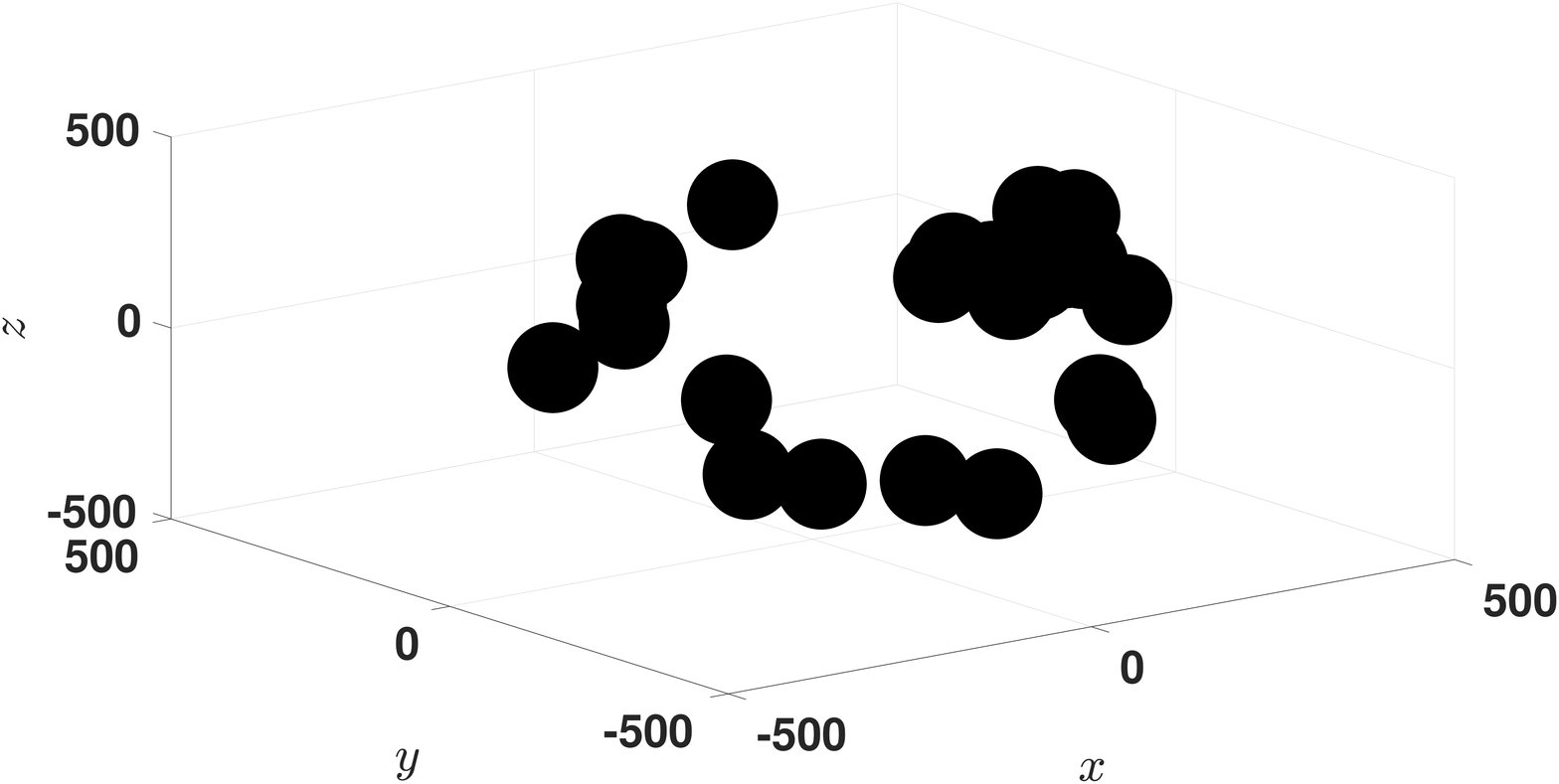}
    \caption{Galaxy with 25 objects after 14000 mega years}
  \end{subfigure}
  %\subfloat[Galaxy with 25 objects after 14000 mega years]{\label{Afig03}\epsfig{figure=figure/GalLF1_final_scatterplot.eps,width=0.48\textwidth}} \quad
  ~\begin{subfigure}[t]{0.4\textwidth}
    \centering
    \includegraphics[width=\textwidth]{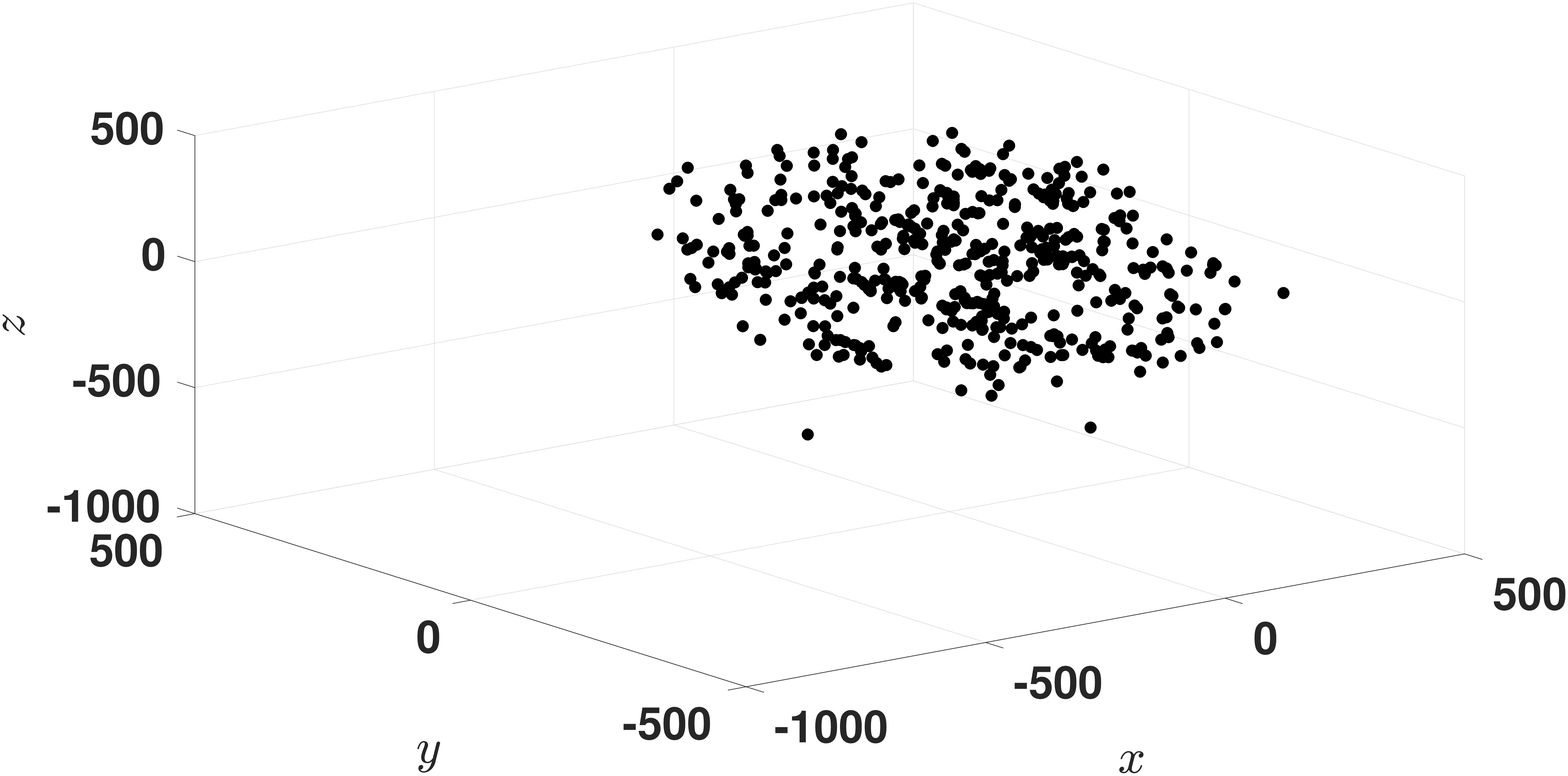}
    \caption{Galaxy with 500 objects after 14000 mega years}
  \end{subfigure}
  %\subfloat[Galaxy with 500 objects after 14000 mega years]{\label{Bfig03}\epsfig{figure=figure/GalHF1_final_scatterplot.eps,width=0.48\textwidth}} \quad
\caption{Galaxy systems with 5 and 500 bodies as the low- and high-fidelity systems; $M=50$; $\Gamma=0$}
\label{fig03}
\end{figure}

\begin{figure}[!tb]
\centering
  \begin{subfigure}[t]{0.4\textwidth}
    \centering
    \includegraphics[width=\textwidth]{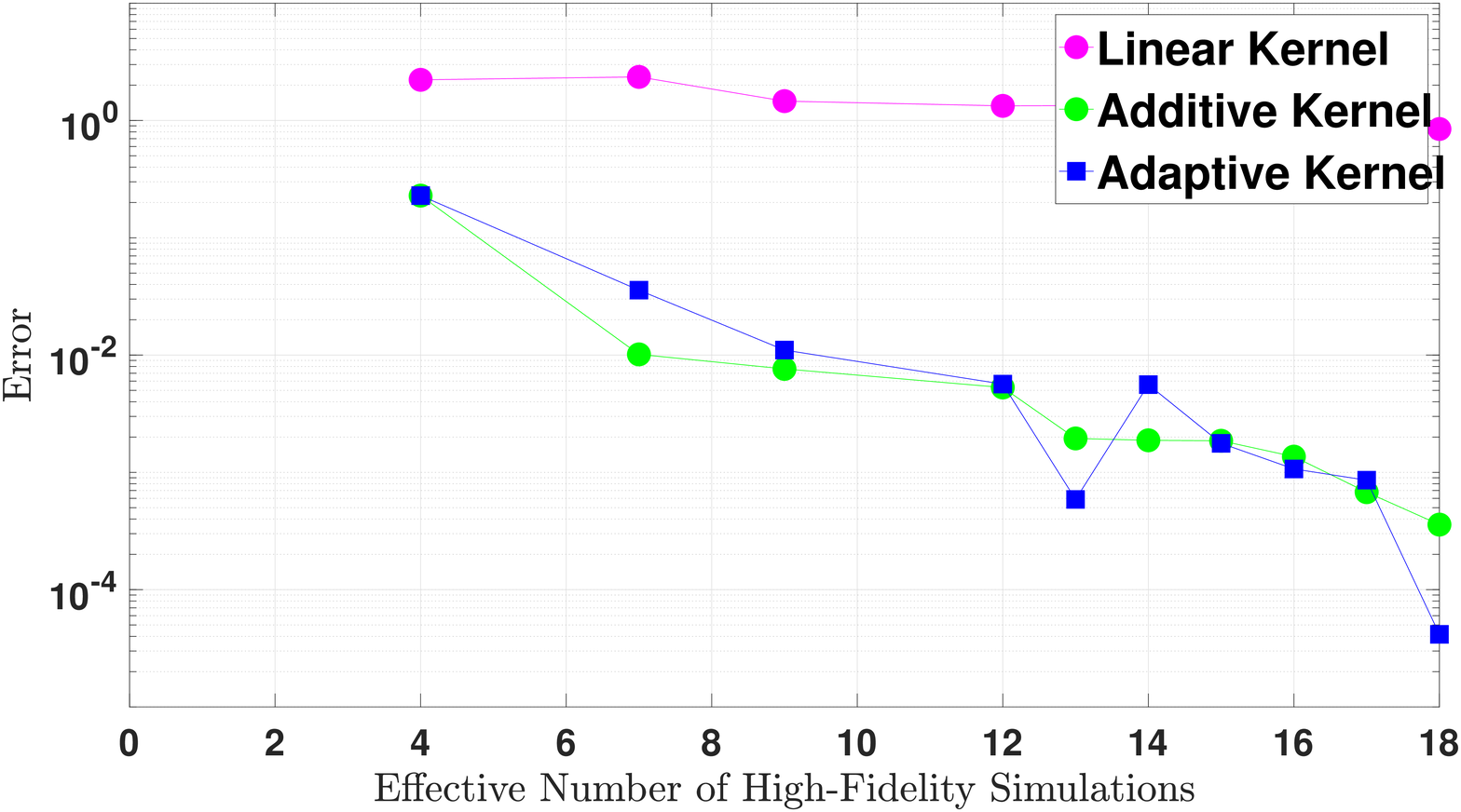}
    \caption{Total Energy}
  \end{subfigure}
  %\subfloat[Total Energy]{\label{Afig3}\epsfig{figure=figure/Gal_TE_corrected.eps,width=0.48\textwidth}} \quad
  ~\begin{subfigure}[t]{0.4\textwidth}
    \centering
    \includegraphics[width=\textwidth]{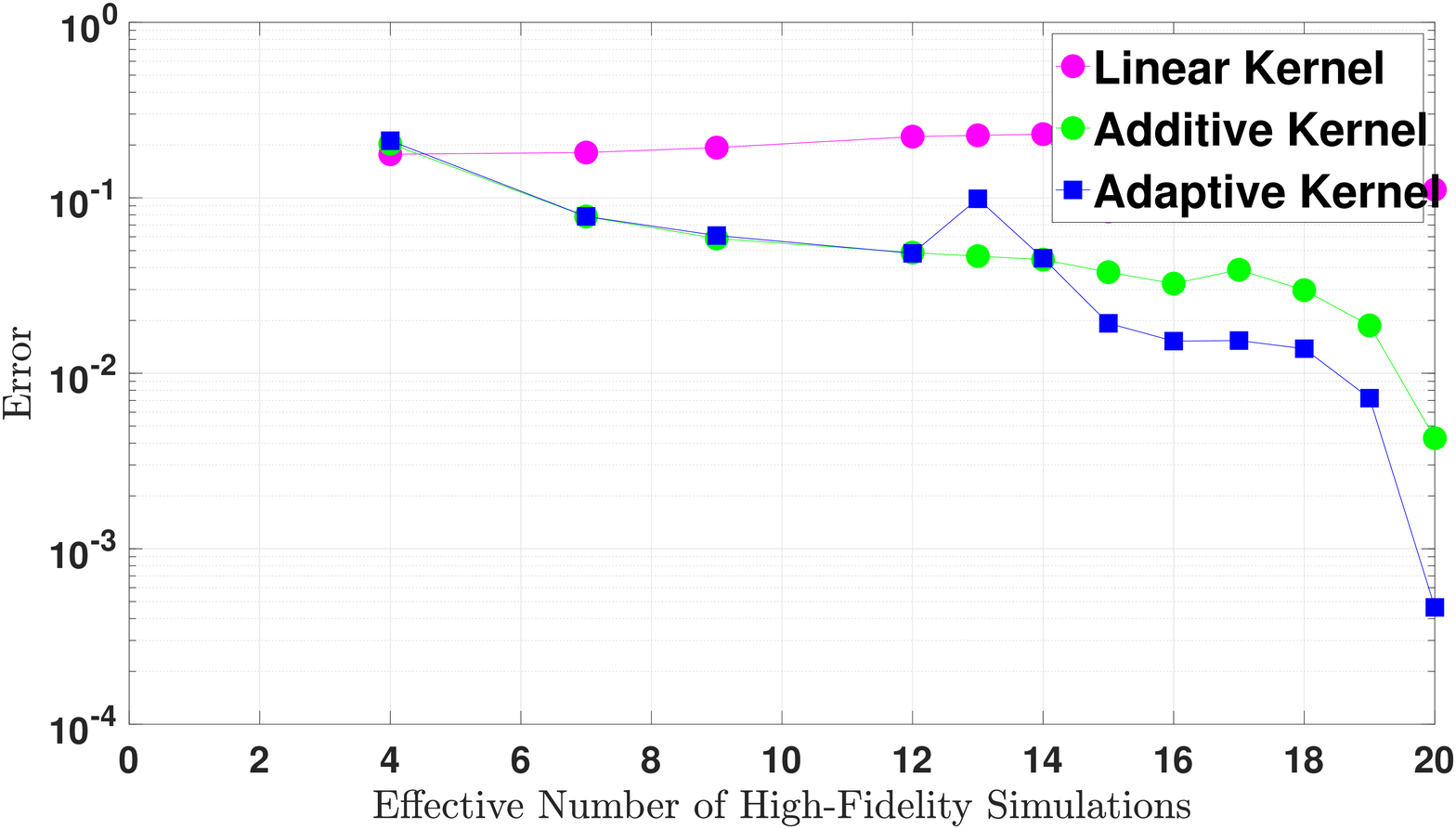}
    \caption{Mean distance in of objects in the $n$-body system relative to the center of the system}
  \end{subfigure}\\
  %\subfloat[mean distance in of objects in the $n$-body system relative to the center of the system]{\label{Bfig3}\epsfig{figure=figure/Gal_R_corrected.eps,width=0.48\textwidth}} \newline\newline\newline\vskip -0.5cm
  \begin{subfigure}[t]{0.4\textwidth}
    \centering
    \includegraphics[width=\textwidth]{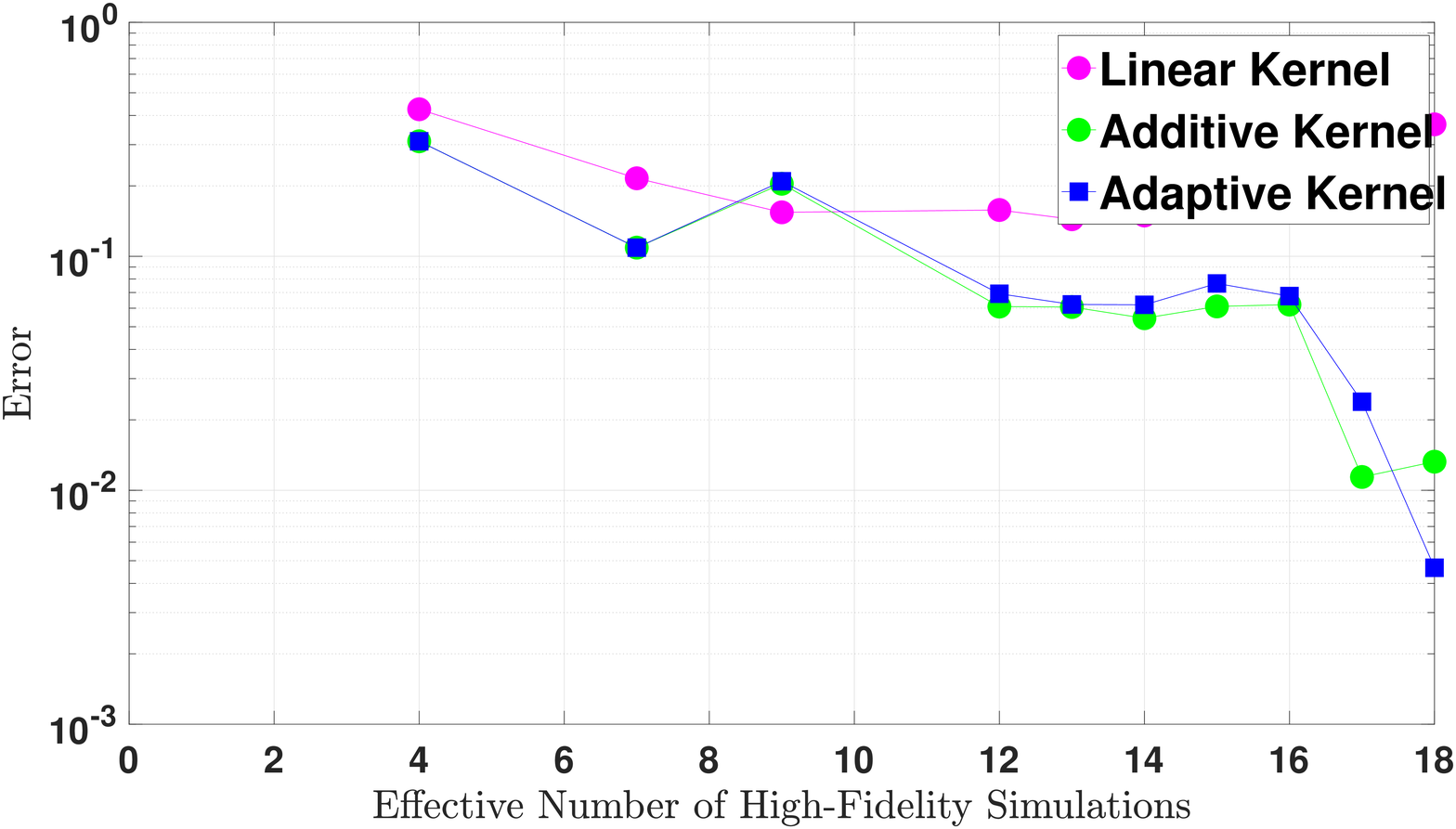}
    \caption{Mean velocity of $n$-body system}
  \end{subfigure}
%\subfloat[mean velocity of $n$-body system]{\label{Cfig3}\epsfig{figure=figure/Gal_R_corrected.eps,width=0.48\textwidth}} \quad
\caption{Median error of the multi-fidelity model constructed based upon the results from the low rank low-fidelity galaxy ($n$ body) model with $n=25$ bodies in the prediction of the properties for the high-fidelity model of the galaxy with $n=500$ bodies after 14000 mega years; Test Problem 2; 36 data points were tested for each case to measure the error.}
\label{fig3}
\end{figure}

The parameter space for this case involve system's total mass ($M\in[50,500]$ solar mass) and the initial rotation ($\Gamma \in [0,0.9]$) about a central axis in the $n$-body cluster imposed on the system. Fidelity here is based on quantized parameters in the model \cite{razi2018quantized}: We consider galaxy systems of 25 and 500 objects as the low- and high-fidelity models, respectively. The quantities of interest for this test problem are total energy, mean distance from the origin, and the mean velocity of the objects in galaxy systems. As illustrated in Fig.~\ref{fig03}, the evolutions of both low- and high-fidelity models is significantly different, and thus we do not expect to be able to use MF simulations to compute detailed trajectories; however, our quantities of interest are averaged metrics, and so using a MF surrogate procedure can be effective. The results provided in Fig.~\ref{fig3} indicate success of of the kernel selection procedures for constructing MF emulators. This level of accuracy is obtained despite the apparent failure of the standard linear kernel MF surrogate. Thus, the kernel selection procedure takes advantages of the optimization process to promote numerical stability and accuracy. Here, the Adaptive Kernel Approach provides us with a slightly more accurate predictive model. This could be the result excluding several ``improper'' choices of kernel function. In the Additive Kernel Approach for this problem, every kernel function in the library has a non-zero contribution in construction of the additive kernel.
  
\revision{
As shown in Fig.~\ref{fig3}, the cost of multi-fidelity simulations scales proportionally to the cost of each added high-fidelity simulation sample. By increasing the number of high-fidelity samples from 4 to 18 for the construction of the multi-fidelity emulator, the prediction error of the emulator models reduces by around two orders of magnitude. This is the case for all outputs of the model, and for two of the proposed approaches. This level of accuracy is achieved despite the complexity of the dynamical system, but is possible by considering only averaged characteristics of the system. We expect significant decline of accuracy if particle trajectories or velocities are directly predicted. 
}

\subsection{Test Problem 3: Associating polymer network} 
Telechelic associated polymer networks in their equilibrium condition can be modeled as a Simulated Gel Network (SGN) with its nodes characterized as ``peers'' and ``superpeers''; for details see \cite{billen2009topological}. The number of links and nodes in the network depends on the temperature. For this test problem, the fractions of superpeers ($n_s$) and peers ($n_p$) at non-dimensional temperature ($T$) of 0.55 based on the Table 1 of Ref. ~\cite{billen2009topological} are set to be 0.604 and 0.396, respectively. The following system of equations and constraints model the network that exhibits the same properties as those of an SGN:
\begin{eqnarray}
n_s\langle k\rangle_s&=& 2l_{ss}+l_{ps},\nonumber\\
n_p \langle k\rangle_s&=& l_{ps}+2l_{pp},\nonumber\\
l_{\text{tot}}&=&l_{pp}+l_{ps}+l_{ss},\nonumber\\
l{ss}&=&(0.85 \pm 0.01)l_{\text{tot}},\nonumber\\
l{pp}&=&(0.02 \pm 0.01)l_{\text{tot}},\nonumber\\
l{ps}&=&(0.13 \pm 0.01)l_{\text{tot}},
\label{eq:NET1}
\end{eqnarray}
where, $l_{ss}$, $l_{ps}$, and $l_{pp}$ represent the probability of having a connection between two superpeers, between one peer and a superpeer, and between two peers, respectively. Here $\langle k\rangle$ denotes the average degree (number of connections for each node), and for $T=0.55$ is set to be $\langle k\rangle =7.131$ and $\langle k\rangle = 0.997$ for superpeers and peers, respectively. We consider 1200 nodes for the model network and randomly generate links based on the calculated values of $l_{ss}$, $l_{ps}$, and $l_{pp}$.

\begin{figure}[!tb]
  \centering
  \begin{subfigure}[t]{0.4\textwidth}
    \centering
    \includegraphics[width=\textwidth]{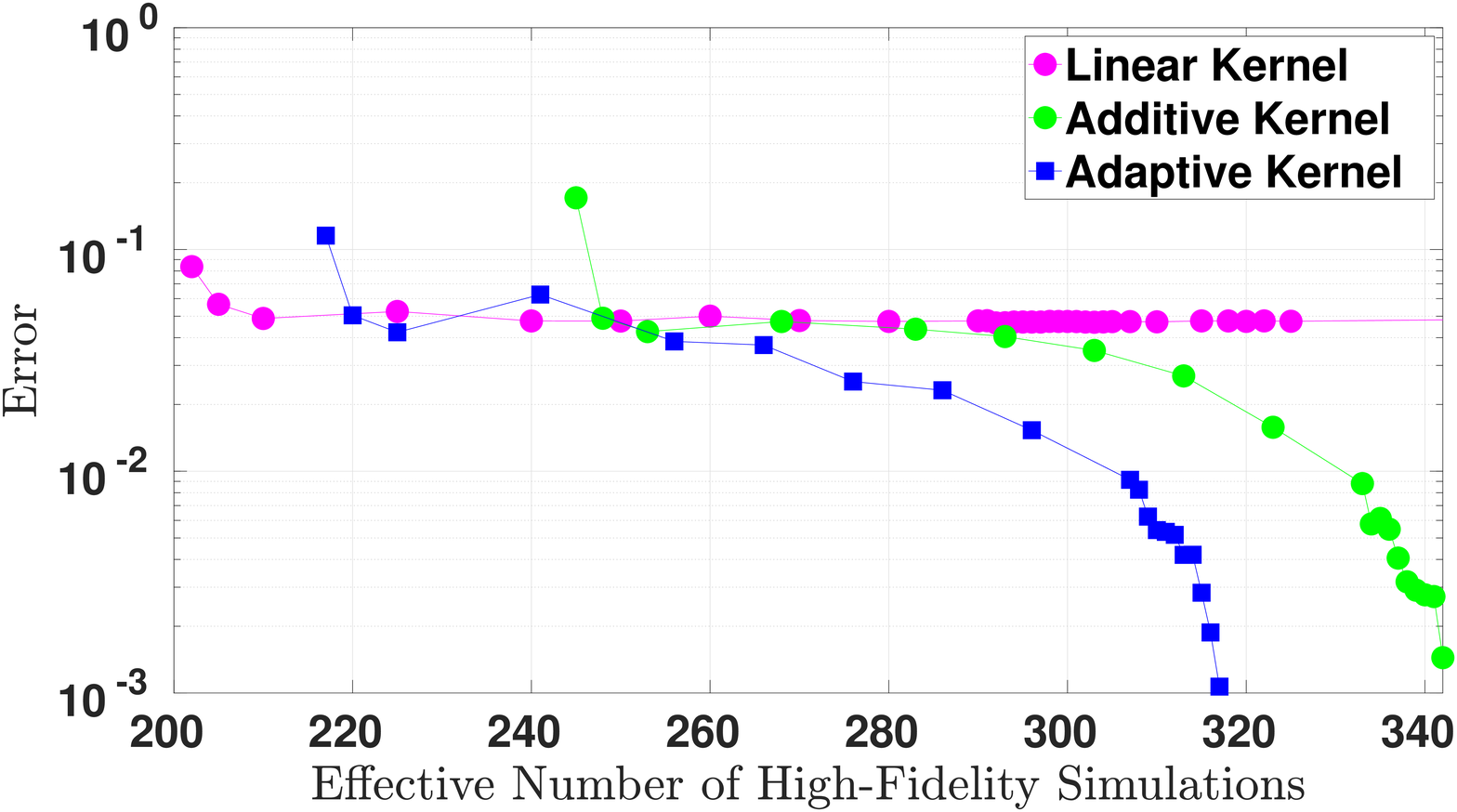}
    \caption{strain field mean}
  \end{subfigure}
  %\subfloat[strain field mean]{\label{Afig7}\epsfig{figure=figure/Pol_Mean_10pop.eps,width=0.48\textwidth}} \quad
  ~\begin{subfigure}[t]{0.4\textwidth}
    \centering
    \includegraphics[width=\textwidth]{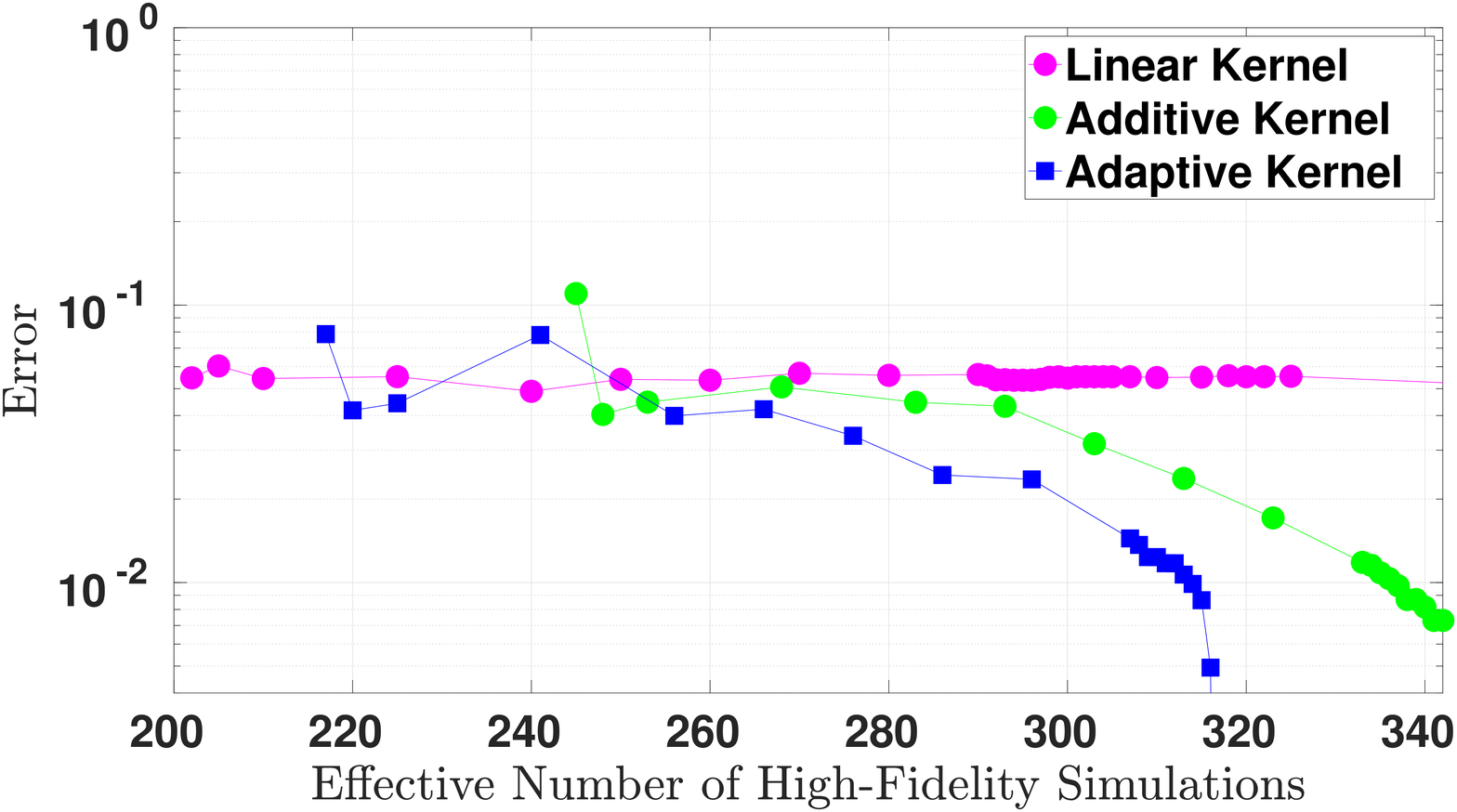}
    \caption{strain field standard deviation}
  \end{subfigure}
  %\subfloat[strain field standard deviation]{\label{Bfig7}\epsfig{figure=figure/Pol_STD_10pop.eps,width=0.48\textwidth}} \quad
\caption{Median error of the multi-fidelity model constructed based upon the results from the low rank low-fidelity model (sparse model polymer network) in the prediction of the mean and standard deviation of the resultant strain field at $t=10$ for the high-fidelity model (dense model polymer network); Test Problem 3; 200 data points were tested for each case to measure the error.}
\label{fig7} 
\end{figure}

In order to study the topological change in the simulated associating polymers under constant stress, we consider the following governing equations for the nodal strain evolution in the network~\cite{christensen2012theory}: 
\begin{equation}
\frac{d\varepsilon}{dt}=\frac{1}{\eta}A\sigma_0,
\label{eq:NET2}
\end{equation}
where $\eta$, $\sigma_0$, $\varepsilon$ and $A$ are a coefficient of viscosity, vectors of constant stress, time-dependent strain and the network's adjacency matrix, respectively. For this test problem, we consider the of effective resistances between any two nodes in a network and use a graph sparsification process~\cite{spielman2011graph} to obtain a low-fidelity network model. Hence, the high- and low-fidelity models in this case are fully dense and sparsified networks, respectively. Here, the parameter space is obtained by uniform random variation of the viscosity coefficient and constant in $[0.05,2.5]$ and between 0 and 2.5, respectively. Two scalar quantities of interest for this test problem are the mean and standard deviation of the resulting strain field at $t=10$. In this case the output quantities of interest are two-dimensional vectors. Once again, as shown in Fig.~\ref{fig7}, the dimension of the outputs being small directly impacts the accuracy of the standard (linear kernel) multi-fidelity model. However, the MF models using the new modified kernel approaches perform much better. Similar to the result of the previous test problems, the adaptive approach appears to be more cost-effective and results in a more accurate predictive model. Once again, this better accuracy is the result of selection of the optimal kernel function exclusively in the adaptive approach instead of estimation of the contribution factor for each kernel function in the function library in the additive approach.

\revision{
Here, as the number of nodes in both low- and high-fidelity networks are the same, the computational cost of the simulation for both are roughly the same. This is the main reason behind the large number of effective high-fidelity simulations required to obtain a reasonable level of accuracy for the MF emulators.  
}

\subsection{Test Problem 4: Plasmonic nano-particle arrays}\label{ssec:plasmonic}
In a wide range of real-world technical problems in engineering, the quantities of interest are scalar quantities. As mentioned previously in this paper, using the linear kernel function in such a situation is essentially ineffective, producing an MF emulator with only 1 degree of freedom, and results in a surrogate model with a significant inaccuracy. Designing a plasmonic nano-particle array with an optimal scattering and/or extinction efficiency is among these real-world scenarios. For the calculation of these efficiency parameters as well as simulation of the optical response of a set of identical non-magnetic metallic nano-spheres with sizes much smaller than the wavelength of light (here 25 nm) the Coupled Dipole Approximation (CDA) method~\cite{Guerin2006JOSAA} is used. Using this method, for $N$ metallic nano-particles described by the same volumetric polarizability $\alpha(\omega)$ and located at vector positions $r_i$, one can compute the local field $E_{loc}(r_i)$ as:
\begin{equation}\label{eqFoldyLax}
\left(1-\frac{\alpha k^2}{\epsilon_0} \sum_{j=1,j\neq i}^{N} \mathbf{\tilde{G}}_{ij}\right) \mathbf E_{loc}(\mathbf r_i)= \mathbf E_0{\mathbf r_i},
\end{equation}
where $\mathbf E_0(r_i)$ is the incident field in vector form, $k$ is the wavenumber in the background medium, $\epsilon_0$ is the dielectric permittivity of vacuum ( $\epsilon_0 = 1$ in CGS unit system), and $\mathbf{\tilde{G}}_{ij}$ is constructed from $3\times3$ block of the overall $3N\times3N$ Green's matrix for the $i$th and $j$th particles and the summation runs through all $j$th particles except for $j=i$.

\begin{figure}[!tb]
  \centering
  \begin{subfigure}[t]{0.4\textwidth}
    \centering
    \includegraphics[width=\textwidth]{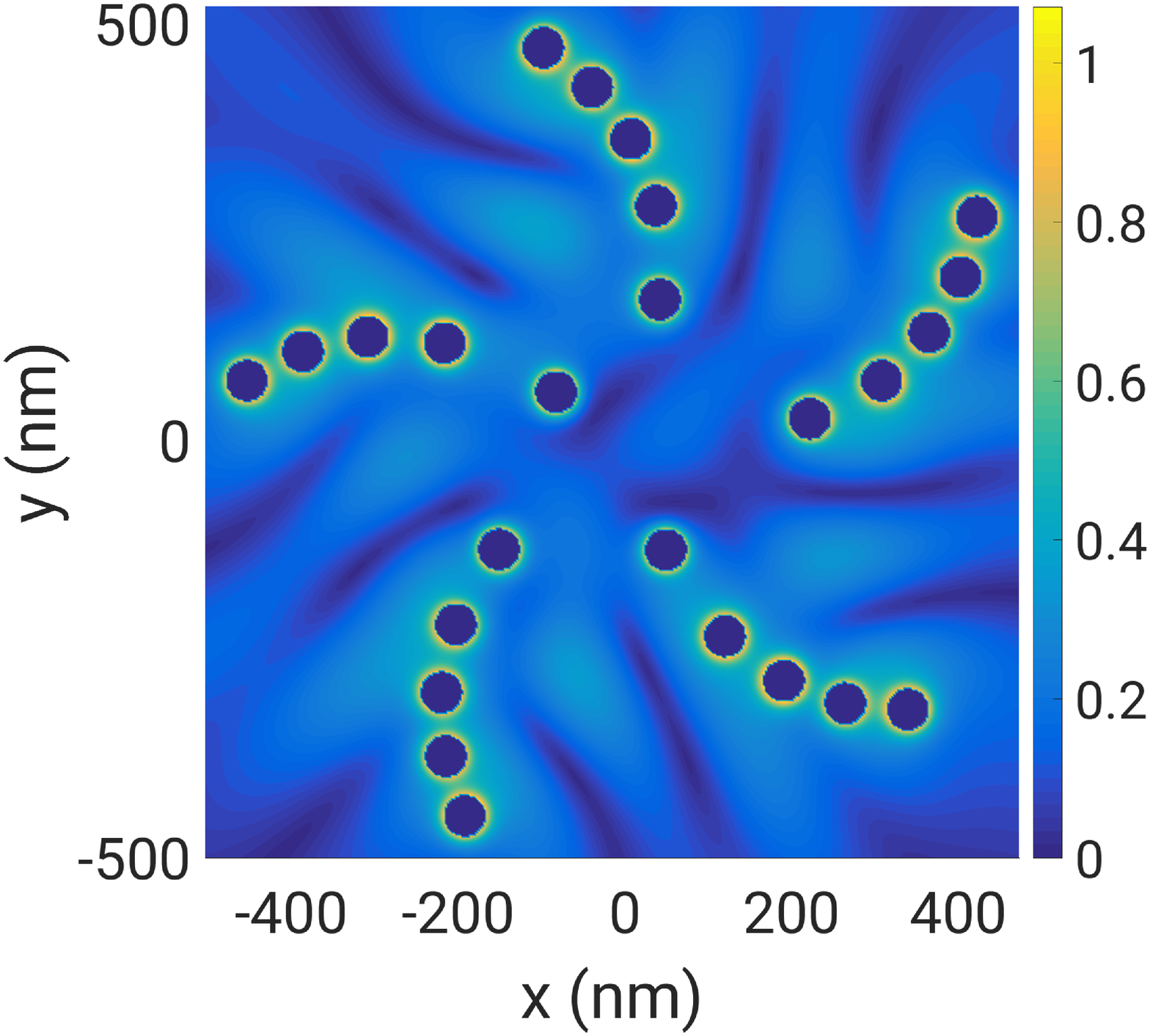}
    \caption{$N=25$}
  \end{subfigure}
  %\subfloat[$N=25$]{\label{Afig9}\epsfig{figure=Plots/Plasmonic_AG_N25.eps,width=0.48\textwidth}} \quad
  ~\begin{subfigure}[t]{0.4\textwidth}
    \centering
    \includegraphics[width=\textwidth]{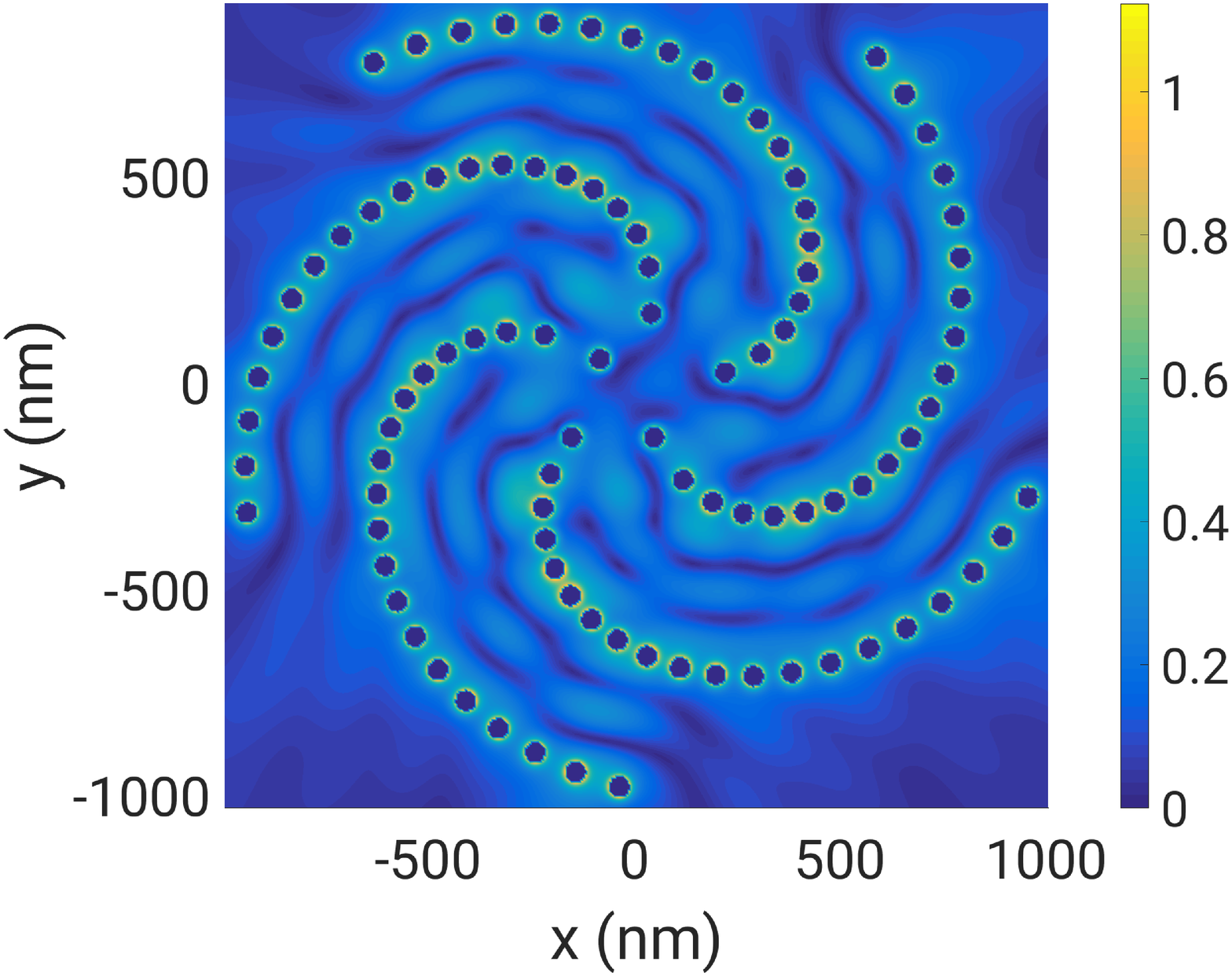}
    \caption{$N=100$}
  \end{subfigure}
  %\subfloat[$N=100$]{\label{Bfig9}\epsfig{figure=Plots/Plasmonic_AG_N100.eps,width=0.48\textwidth}}
\caption{Calculated local extinction field from CDA; Test Problem 4; silver; $\lambda=513.228$ nm, $a_{vs}=100.019$ and $\alpha_{vs}=2.5358$ rad.}
\label{fig9}
\end{figure}

\begin{figure}[!tb]
\centering
  \begin{subfigure}[t]{0.4\textwidth}
    \centering
    \includegraphics[width=\textwidth]{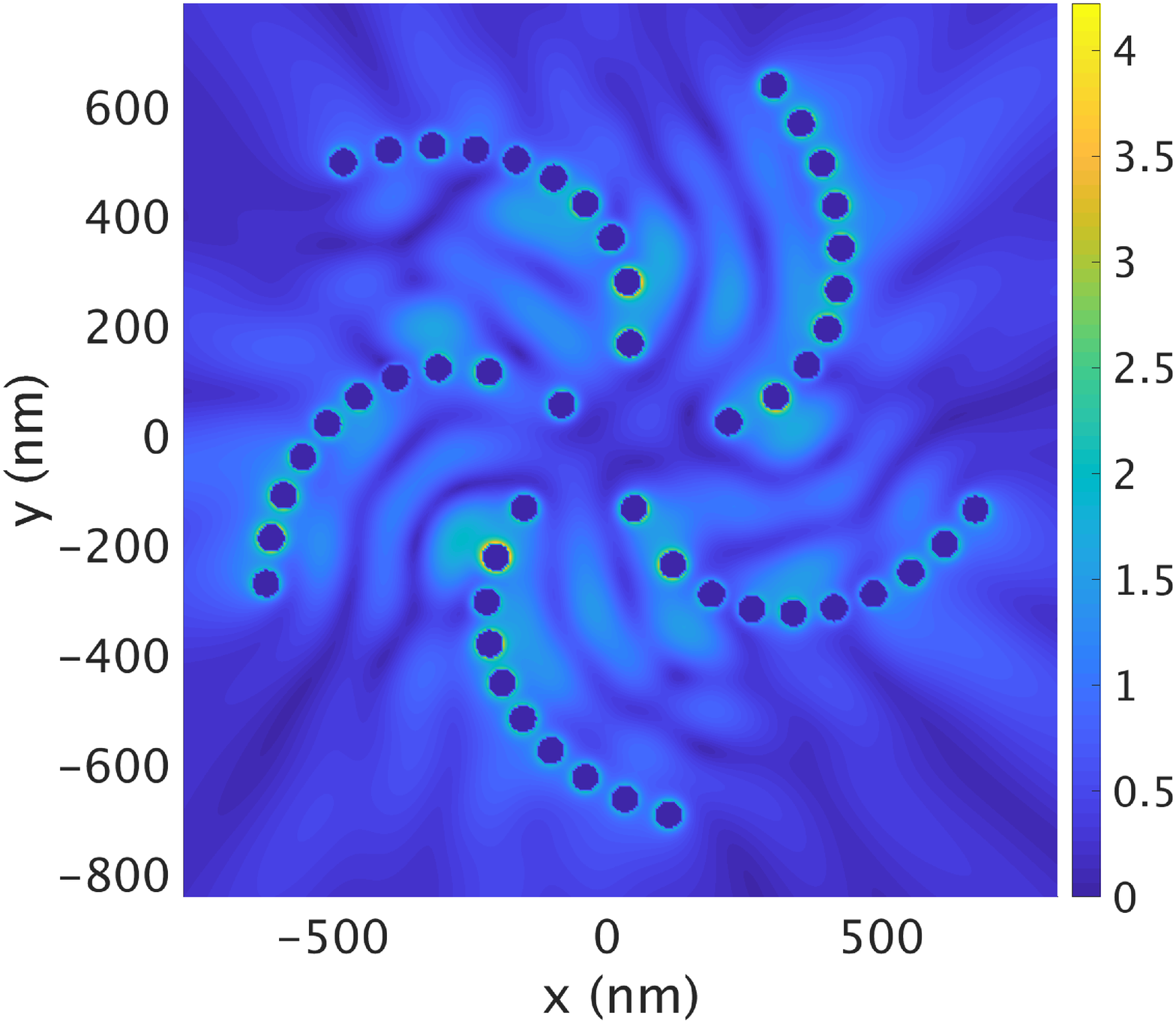}
    \caption{$N=50$}
  \end{subfigure}
  %\subfloat[$N=50$]{\label{Afig99}\epsfig{figure=figure/Al_Vogel_N50.eps,width=0.48\textwidth}} \quad
  ~\begin{subfigure}[t]{0.4\textwidth}
    \centering
    \includegraphics[width=\textwidth]{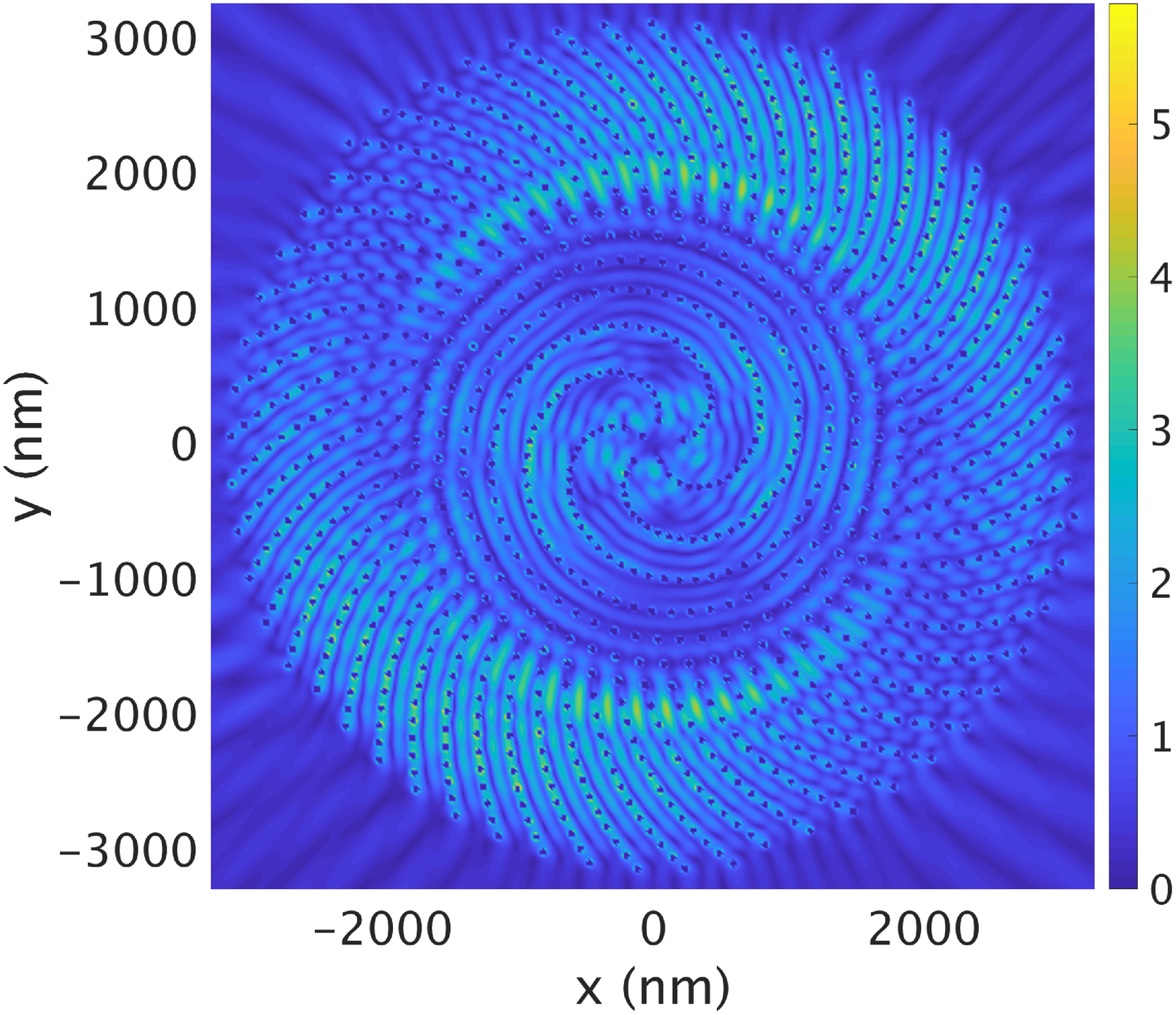}
    \caption{$N=1000$}
  \end{subfigure}
  %\subfloat[$N=1000$]{\label{Bfig99}\epsfig{figure=figure/Al_Vogel_N1000.eps,width=0.48\textwidth}}
\caption{Calculated local extinction field from CDA; Test Problem 4; Aluminum; $\lambda=213.228$ nm, $a_{vs}=100.019$ and $\alpha_{vs}=2.5358$ rad.}
\label{fig99}
\end{figure}

\begin{figure}[!tb]
  \centering
  \begin{subfigure}[t]{0.4\textwidth}
    \centering
    \includegraphics[width=\textwidth]{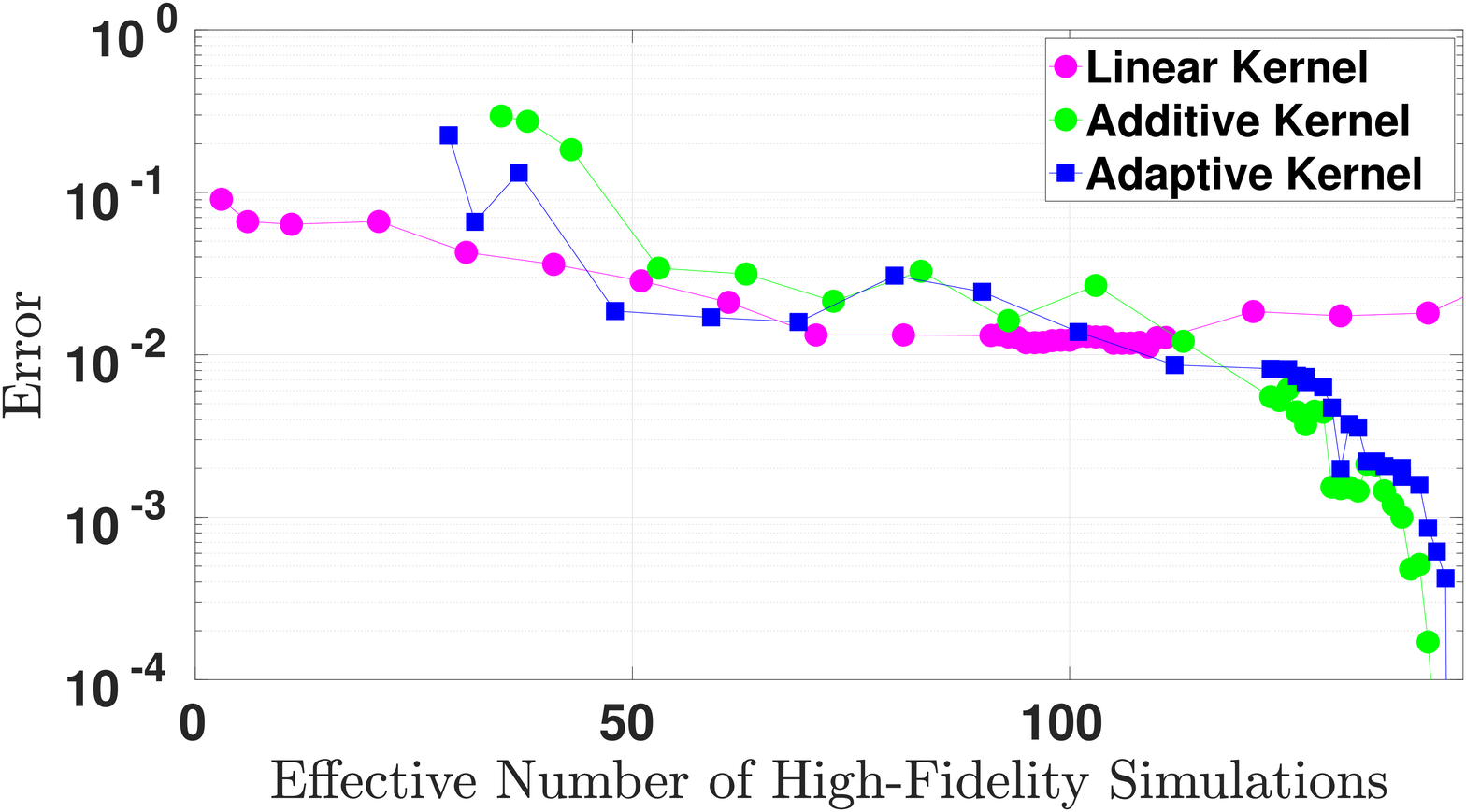}
    \caption{$Q_{ext}$}
  \end{subfigure}
  %\subfloat[$Q_{ext}$]{\label{Afig4}\epsfig{figure=figure/plasmonic_Qext_Ag_50vs1000_250samples.eps,width=0.48\textwidth}} \quad
  ~\begin{subfigure}[t]{0.4\textwidth}
    \centering
    \includegraphics[width=\textwidth]{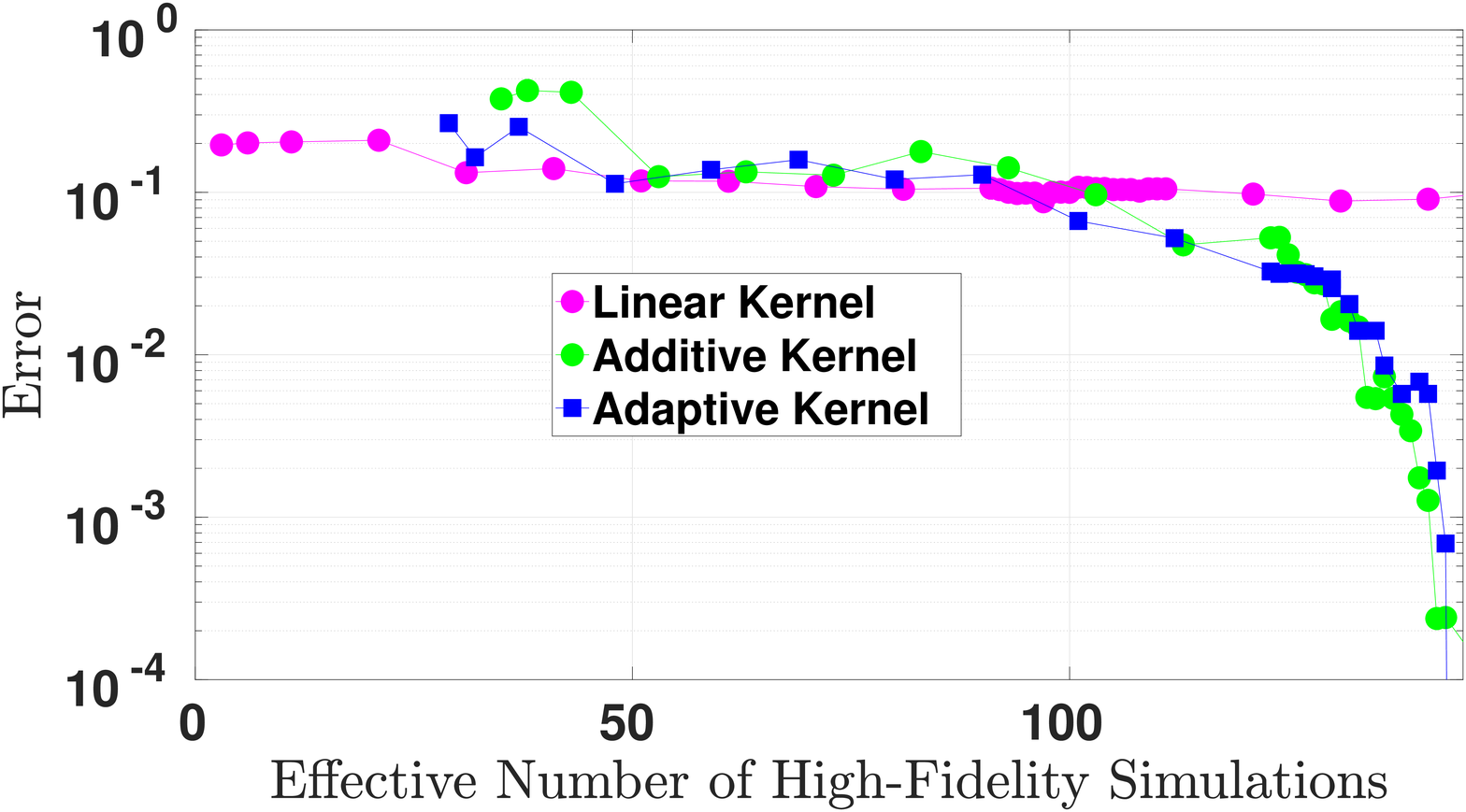}
    \caption{$Q_{sc}$}
  \end{subfigure}
  %\subfloat[$Q_{sc}$]{\label{Bfig4}\epsfig{figure=figure/plasmonic_Qsc_Ag_50vs1000_250samples.eps,width=0.48\textwidth}} \quad
\caption{Median error of the multi-fidelity model constructed based upon the results from the low rank low-fidelity model with 25 nano-particles ($N_L=25$) in the prediction of the extinction and scattering efficiencies for the high-fidelity models with 100 silver nano-particles ($N_H=100$); Test Problem 4; 250 data points were tested for each case to measure the error.}
\label{fig4}
\end{figure}

\begin{figure}[!tb]
\centering
  \begin{subfigure}[t]{0.4\textwidth}
    \centering
    \includegraphics[width=\textwidth]{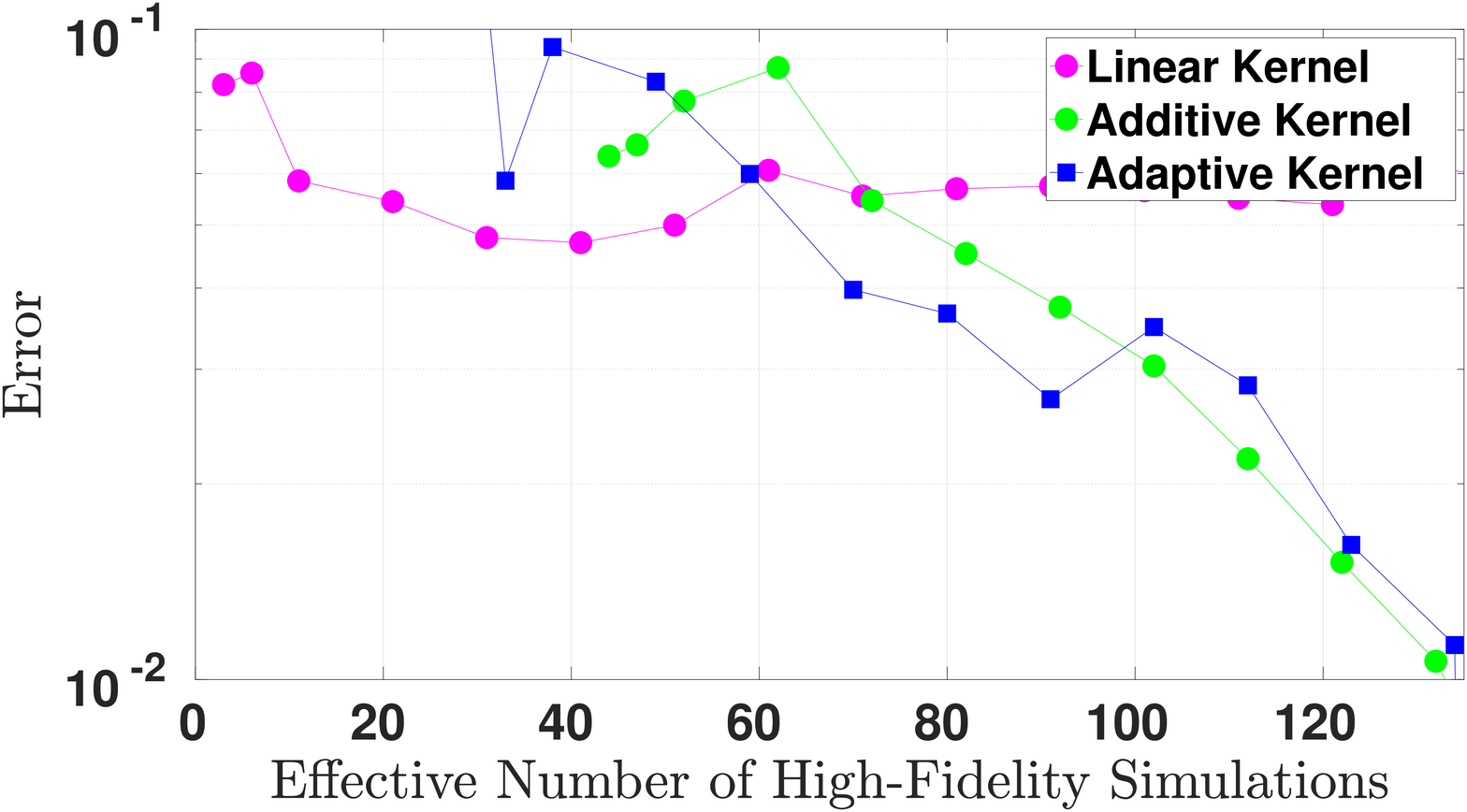}
    \caption{$Q_{ext}$}
  \end{subfigure}
  %\subfloat[$Q_{ext}$]{\label{Afig5}\epsfig{figure=figure/plasmonic_Qext_Al_50vs1000_250samples.eps,width=0.48\textwidth}} \quad
  ~\begin{subfigure}[t]{0.4\textwidth}
    \centering
    \includegraphics[width=\textwidth]{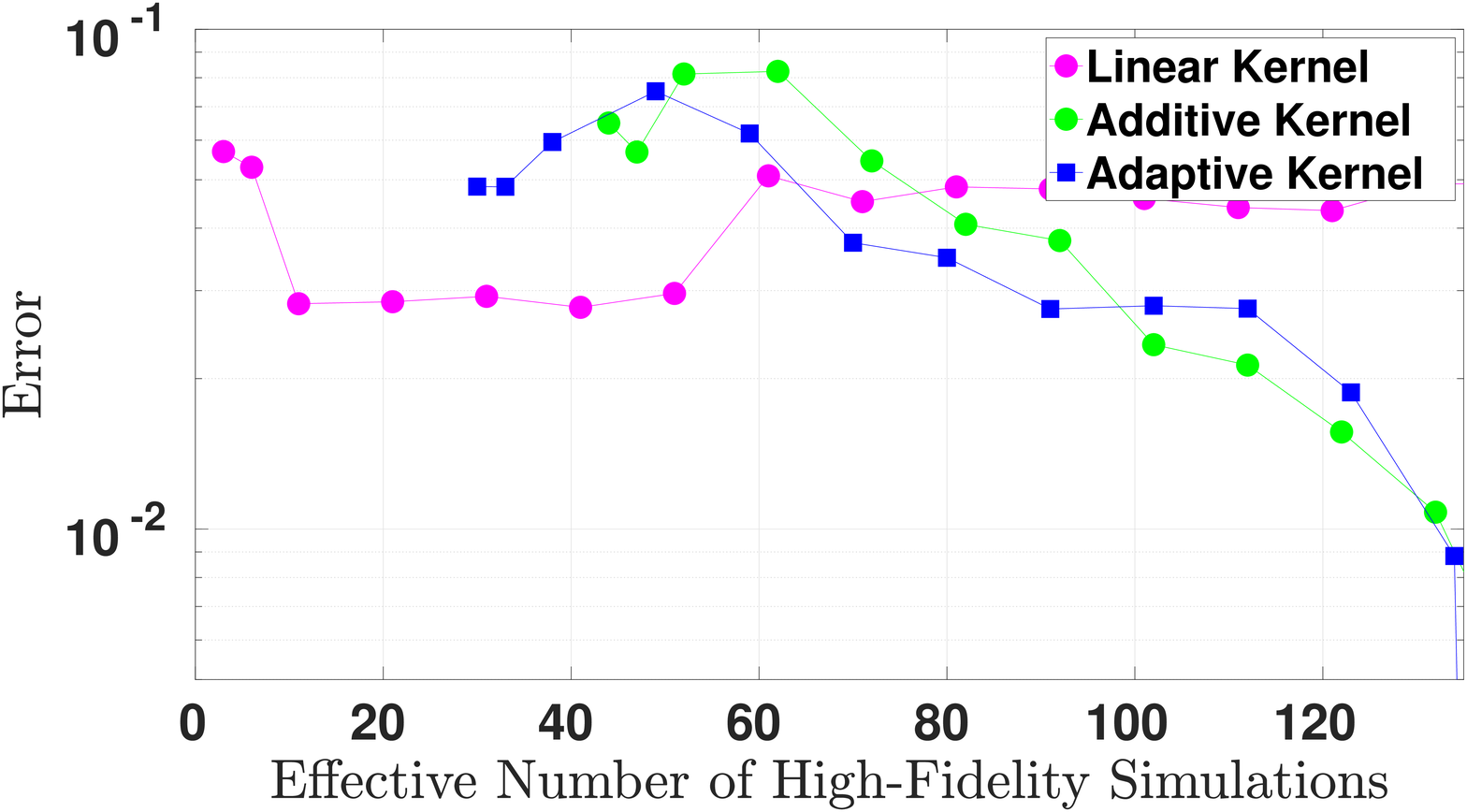}
    \caption{$Q_{sc}$}
  \end{subfigure}
  %\subfloat[$Q_{sc}$]{\label{Bfig5}\epsfig{figure=figure/plasmonic_Qsc_Al_50vs1000_250samples.eps,width=0.48\textwidth}} \quad
\caption{Median error of the multi-fidelity model constructed based upon the results from the low rank low-fidelity model with 50 nano-particles ($N_L=50$) in the prediction of the extinction and scattering efficiencies for the high-fidelity models with 1000 aluminum nano-particles ($N_H=1000$); Test Problem 4; 250 data points were tested for each case to measure the error.}
\label{fig5}
\end{figure}

By solving \eqref{eqFoldyLax}, the total local fields ($\mathbf E_{loc}(r_i)$) are computed, and consequently the local electric dipole polarizations across the array can be determined. The scattering and extinction cross-sections can be computed using the resulting polarizations. Finally, the corresponding extinction ($Q_{ext}$) and scattering ($Q_{sc}$) efficiencies are obtained by normalization of scattering and extinction cross-sections with respect to the total projected area of the array (i.e., sum of the areas of the particles projected perpendicularly to the direction of the excitation beam). 
   
In particular we consider Vogel spirals arrays \cite{Christofi2016OL,razi2018plasmonics}, which can uniquely identify useful nano-particle array configurations with only four parameters: the number of particles, ($\alpha_{vs}$), divergence angle, incident wavelength, and scaling factor. 

Here, the number of particles in the array defines the level of fidelity. Varying this fidelity parameter can result in significantly different outputs (compare the two panels in each of Figs.~\ref{fig9} and Figs.~\ref{fig99}). Also the limited dimension of the output parameter space for this case makes the choice of the linear kernel function prohibitive. The results shown in Figs.~\ref{fig4} and ~\ref{fig5} confirm this deduction for MF models of silver and aluminum nano-particle arrays. The number $N$ of particles in the low- and high-fidelity models for both cases are $N = 25$ and $N = 100$, respectively for Figure \ref{fig9}, and are $N = 50$ and $N = 1000$ for Figure \ref{fig99}, respectively. Both of the novel kernel selection approaches are successfully applied to develop surrogate models with much higher accuracy (see Figs.~\ref{fig4} and ~\ref{fig5}). The computational cost of the kernel selection optimization for this test problem are insignificant in comparison to the cost of running one high-fidelity simulation. Hence, the prediction error of the multi-fidelity emulator coverages more smoothly and monotonically. Here, the accuracy of both new kernel selection MF methods appear to be similar.
  
\revision{
Due to the complexity of plasmonic models for large-scale arrays of particles, the prediction accuracy is lower in Figs.~\ref{fig9} compared to Figs.~\ref{fig99}. Another reason for this lower accuracy is the more pronounced difference between parametric dependence of the low- and high-fidelity nanoparticle arrays and corresponding plasmonic models. Similar parametric dependence between low- and high-fidelity models is beneficial for the selection of a better set of high-fidelity sampling points. 
}

\subsection{Test Problem 5:  An incompressible flow in two-dimensional channels with stenosis}\label{ssec:stenosis}
\newcommand*{\everymodeprime}{\ensuremath{\prime}}
Our previous examples have investigated the two proposed kernel selection schemes when the model output space has small dimension. In most of the cases above, the adaptive kernel approach has a small advantage. However, it is also necessary to investigate the performance of these approaches for cases in which the dimension of output space is high enough to pick sufficient number of ``important'' samples for the construction of the MF model even with the linear kernel function. The objective of considering this test problem is to provide some insights about the competitiveness of the adaptive kernel approach compared to the linear kernel function method, when the Gramian matrix for both is not ill-conditioned. \rrevision{Since the previous examples suggest that the adaptive kernel procedure is slightly more accurate than the additive method, we focus on highlighting results for the adaptive kernel only. Note that the additive kernel is slightly less computationally expensive than the adaptive kernel method.}
   
\begin{figure}[!tb]
  \centering
  \begin{subfigure}[t]{0.4\textwidth}
    \centering
    \includegraphics[width=\textwidth]{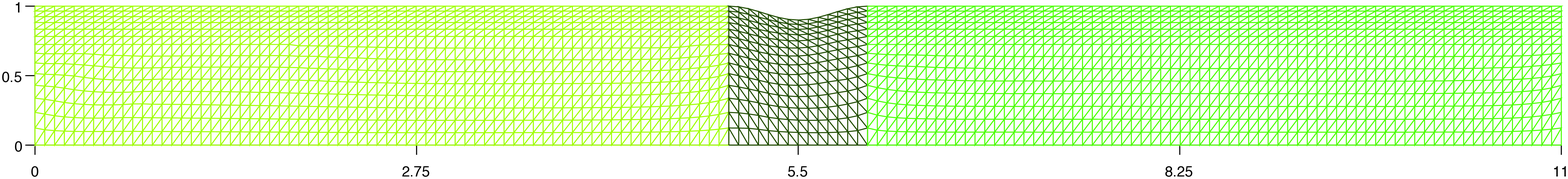}
    \caption{$f_0=0.1$}
  \end{subfigure}
%\subfloat[$f_0=0.1$]{\label{Cfig66}\epsfig{figure=figure/f01.eps,width=0.9\textwidth}} \quad
  ~\begin{subfigure}[t]{0.4\textwidth}
    \centering
    \includegraphics[width=\textwidth]{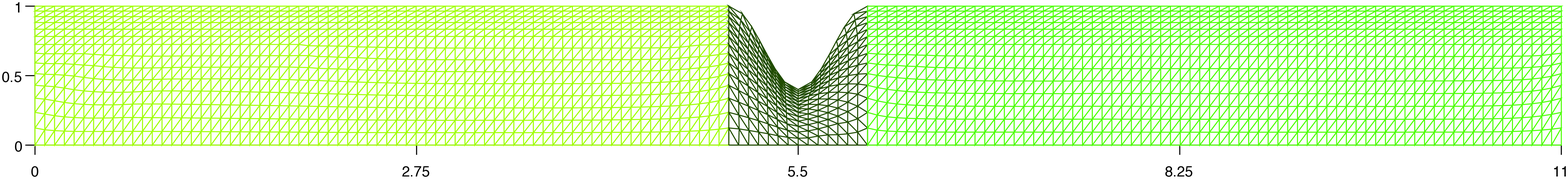}
    \caption{$f_0=0.6$}
  \end{subfigure}
%\subfloat[$f_0=0.6$]{\label{Dfig66}\epsfig{figure=figure/f06.eps,width=0.9\textwidth}} \quad
\caption{Different geometrical configurations tested for two-dimensional stenotic flows; stenosis shape function dependence on the axial coordinate $x$: $f(x)=(1-\frac{f_0}{2}\left(1+\cos{\left(2\pi(x-5)\right)}\right)$; Test Problem 5}
\label{fig66}
\end{figure}

In many fluid dynamics problems a large amount of data can be extracted from a solution, and thus are a good example of high-dimensional output space. Here, we consider a two-dimensional incompressible stenotic flow in a large channel with a length and width of 11$m$ and 1 $m$, respectively. The shape of the curve of stenosis in the middle of the channel follows the equation (see Fig~\ref{fig66}):
\begin{equation}
y(x,y_0)=1-\frac{y_0}{2}\left(1+\cos{\left(2\pi\left(x-5\right)\right)}\right).
\label{eq:P5_1}
\end{equation} 
The parameter space for this problem is two-dimensional, comprised of the flow Reynolds number (Re) and the parameter $y_0$: $(\text{Re},y_0) \in [0.1,0.6]\times[10,500]$. The inlet flow stream has a velocity of 1$m/s$ and for simplicity the fluid kinematic viscosity $\nu$ is set to be $1/ \text{Re}$. The goal of constructing predictive models for this problem is to estimate horizontal and vertical velocity ($u$ and $v$, respectively) profiles in the middle of the channel ($x=5.5$m) as well as the shear stress on top and bottom wall of the channel ($y=0$ and 1, respectively) at $t=30$ sec. The HF model in this case solves the full unsteady incompressible Navier-Stokes equation, which can be written in vector form as follows~\cite{anderson2016computational}:
\begin{equation}
\frac{\partial U}{\partial t}+\frac{\partial E}{\partial x} + \frac{\partial F}{\partial y}  + H= 0,
\label{eq:P5_2}
\end{equation} 
where
\begin{eqnarray}
U&=&[0,u,v]^T,\nonumber\\
E&=&[u,u^2+P - \nu \frac{\partial u}{\partial x},uv-\nu \frac{\partial v}{\partial x}]^T,\nonumber\\
F&=&[v,uv- \nu \frac{\partial u}{\partial y},v^2+P-\nu \frac{\partial v}{\partial y}]^T,\nonumber\\
H&=&[-f_x,-f_y,-f_z]^T\nonumber.
\label{eq:P5_3}
\end{eqnarray}
Here, $P$ and $f$ denote pressure and external force, respectively. The LF model is the coupled linearized Navier-Stokes equation, which are obtained by decomposing the solution into a basic state and perturbed state:
\begin{equation}
\frac{\partial V}{\partial t} +W\cdot \nabla{V}+V\cdot \nabla{W}=-\nabla{P}+\nu \nabla^2{V}+f,
\label{eq:P5_4}
\end{equation} 
where $W$ and $W+\varepsilon V$ are the basic and perturbed states, respectively and $\varepsilon<<1$. 
  
\begin{figure}[!tb]
  \centering
  \begin{subfigure}[t]{0.4\textwidth}
    \centering
    \includegraphics[width=\textwidth]{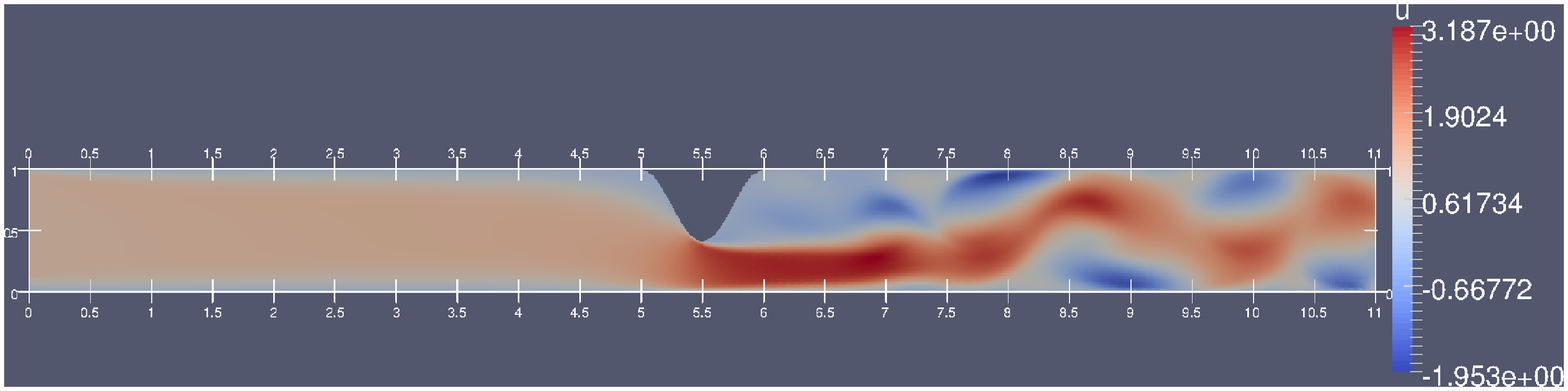}
    \caption{$u$; full Navier-stokes equation}
  \end{subfigure}
%\subfloat[$u$; full Navier-stokes equation]{\label{Afig666}\epsfig{figure=figure/stenosisHF30_u.eps,width=0.9\textwidth}} \quad
  ~\begin{subfigure}[t]{0.4\textwidth}
    \centering
    \includegraphics[width=\textwidth]{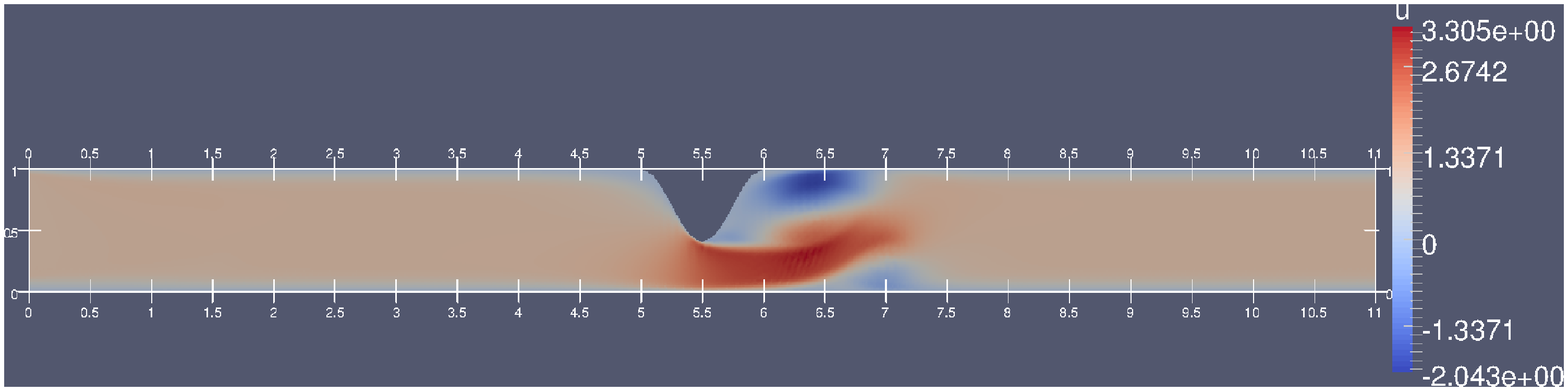}
    \caption{$u$; coupled linearized Navier stokes}
  \end{subfigure}\\
%\subfloat[$u$; coupled linearized Navier stokes]{\label{Bfig666}\epsfig{figure=figure/stenosisLF30_u.eps,width=0.9\textwidth}} \quad
  \begin{subfigure}[t]{0.4\textwidth}
    \centering
    \includegraphics[width=\textwidth]{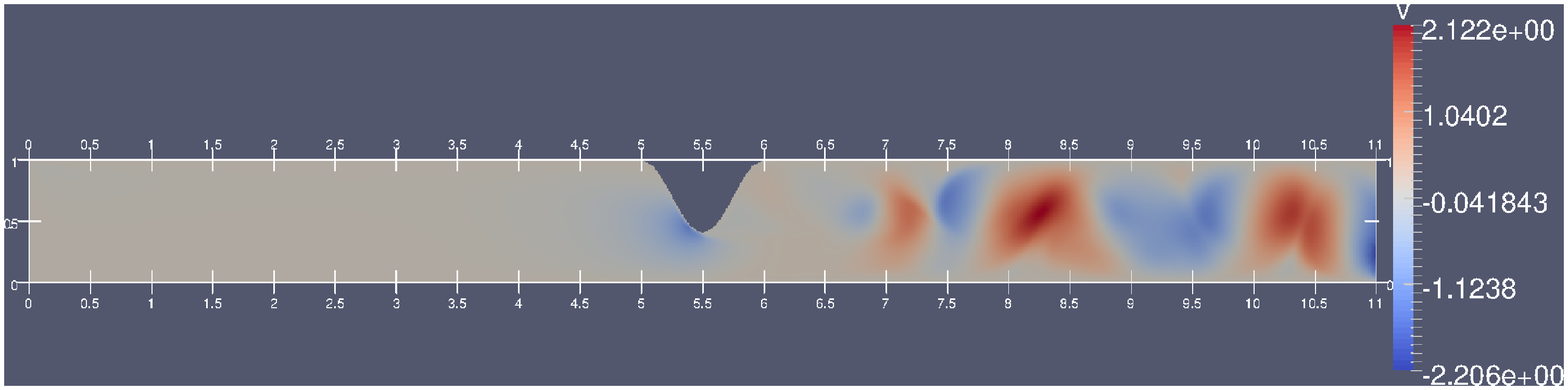}
    \caption{$v$; full Navier-stokes equation}
  \end{subfigure}
%\subfloat[$v$; full Navier-stokes equation]{\label{Cfig666}\epsfig{figure=figure/stenosisHF30_v.eps,width=0.9\textwidth}} \quad
  ~\begin{subfigure}[t]{0.4\textwidth}
    \centering
    \includegraphics[width=\textwidth]{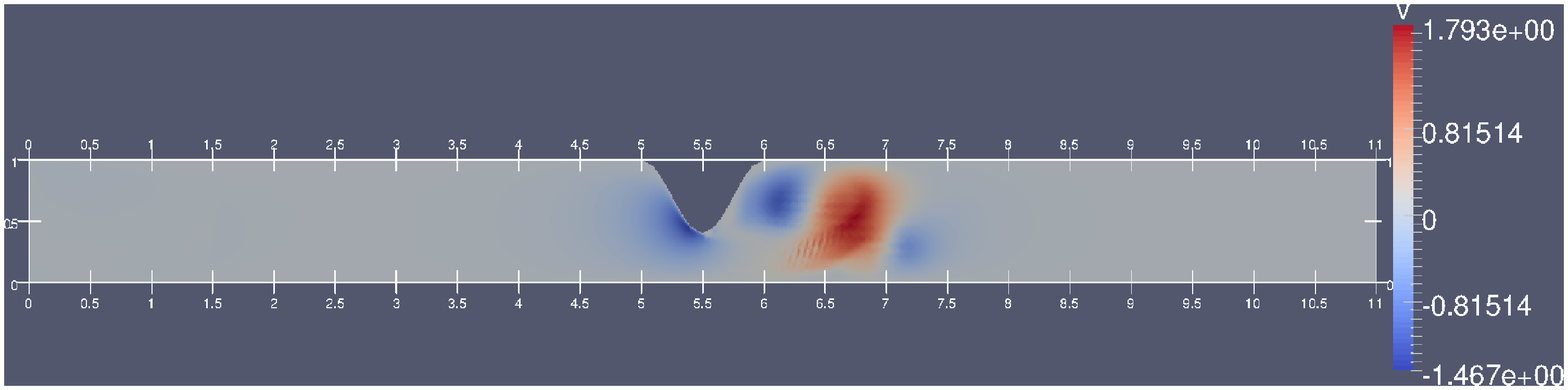}
    \caption{$v$; coupled linearized Navier-stokes}
  \end{subfigure}
%\subfloat[$v$; coupled linearized Navier stokes]{\label{Dfig666}\epsfig{figure=figure/stenosisLF30_v.eps,width=0.9\textwidth}} \quad
\caption{Two-dimensional stenotic flows; $t=30$ sec, $\text{Re}=500$ and $f_0=0.6$; Test Problem 5}
\label{fig666}
\end{figure}

\begin{figure}[!tb]
\centering
  \begin{subfigure}[t]{0.4\textwidth}
    \centering
    \includegraphics[width=\textwidth]{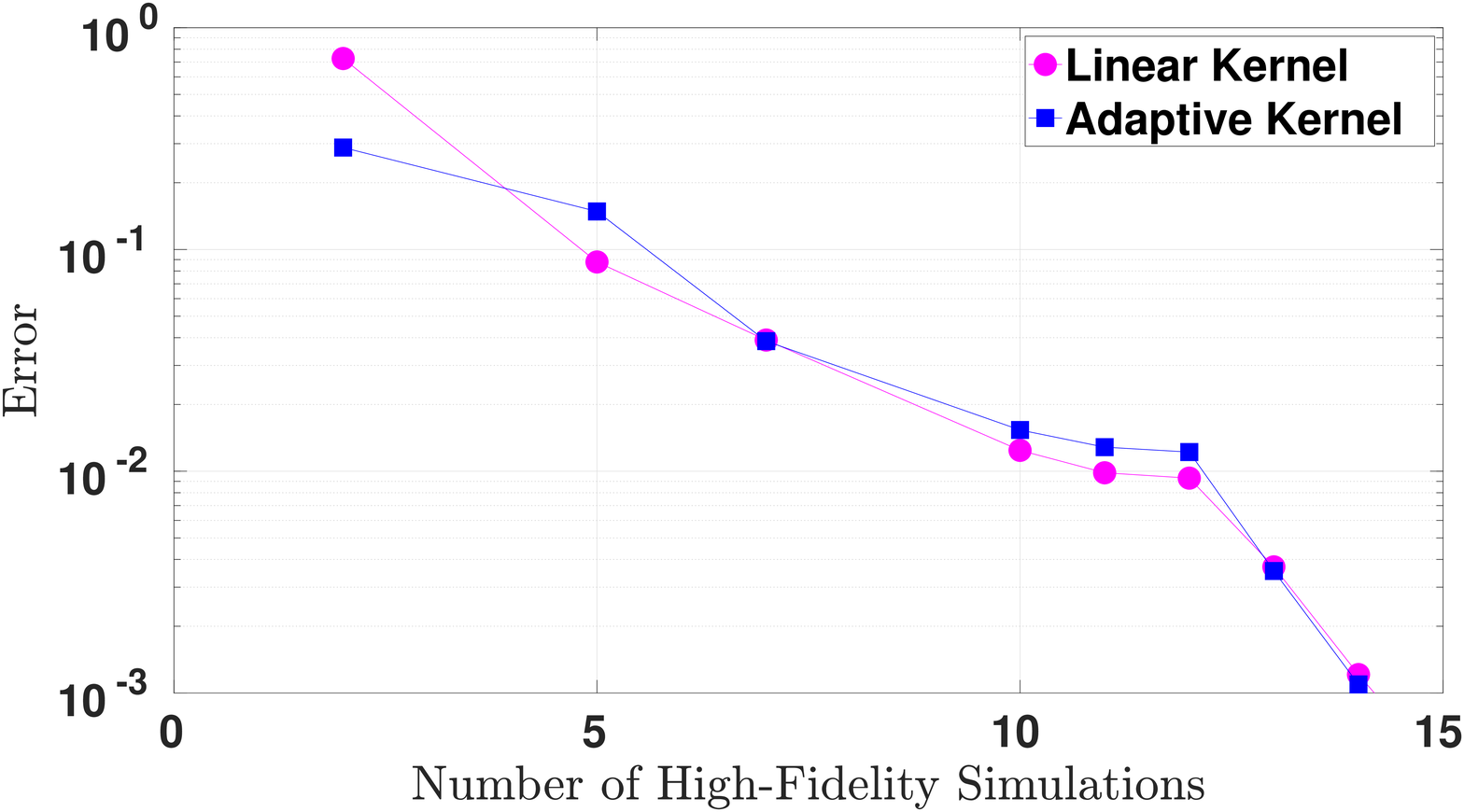}
    \caption{horizontal velocity ($u$) profile at $x=6.5$ m}
  \end{subfigure}
%\subfloat[horizontal velocity ($u$) profile at $x=6.5$ m]{\label{Afig6}\epsfig{figure=figure/Stenosis_u.eps,width=0.48\textwidth}} \quad
  ~\begin{subfigure}[t]{0.4\textwidth}
    \centering
    \includegraphics[width=\textwidth]{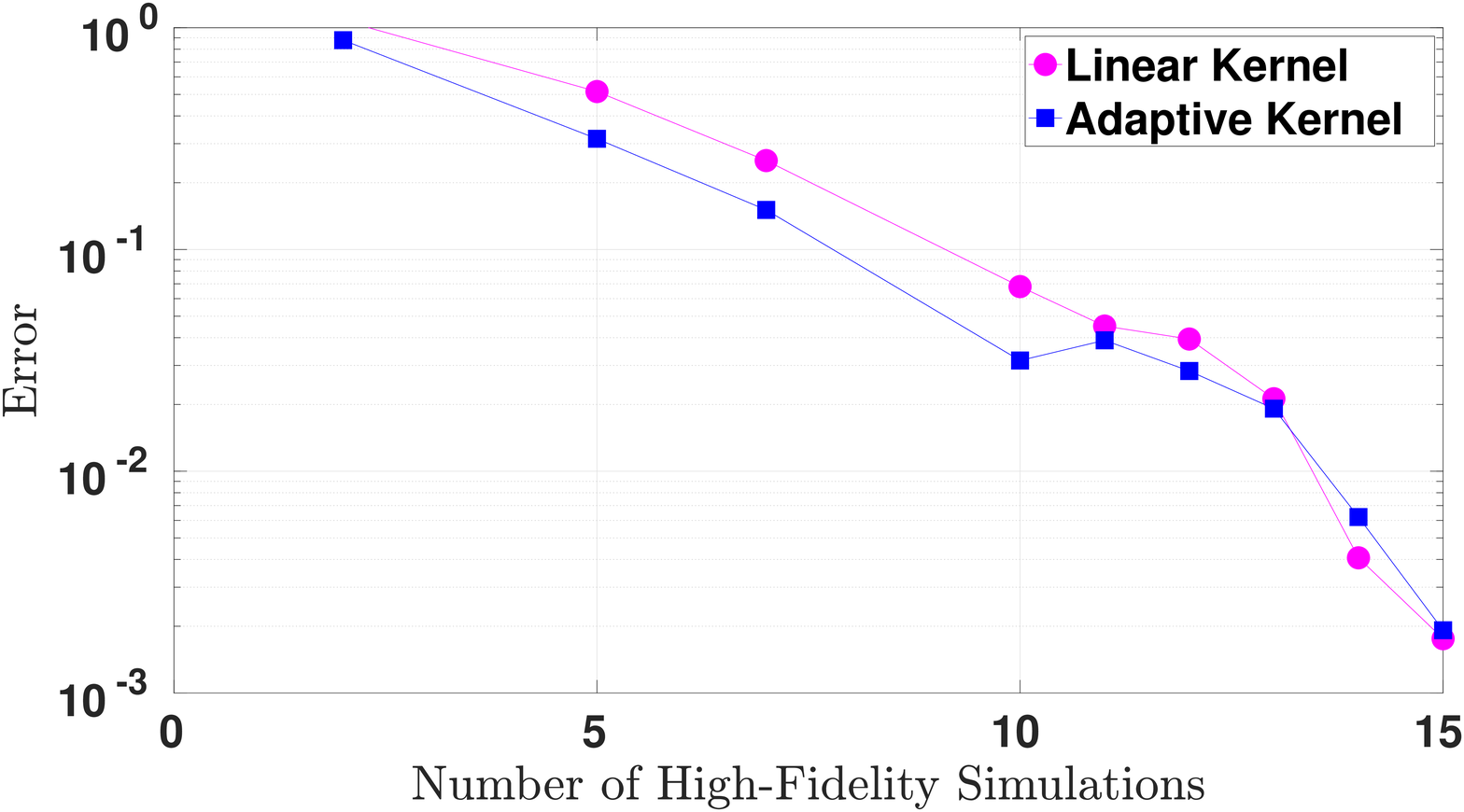}
    \caption{vertical velocity ($v$) profile at $x=6.5$ m}
  \end{subfigure}\\
%\subfloat[vertical velocity ($v$) profile at $x=6.5$ m]{\label{Bfig6}\epsfig{figure=figure/Stenosis_v.eps,width=0.48\textwidth}} \newline\newline\newline\vskip -0.5cm
  \begin{subfigure}[t]{0.4\textwidth}
    \centering
    \includegraphics[width=\textwidth]{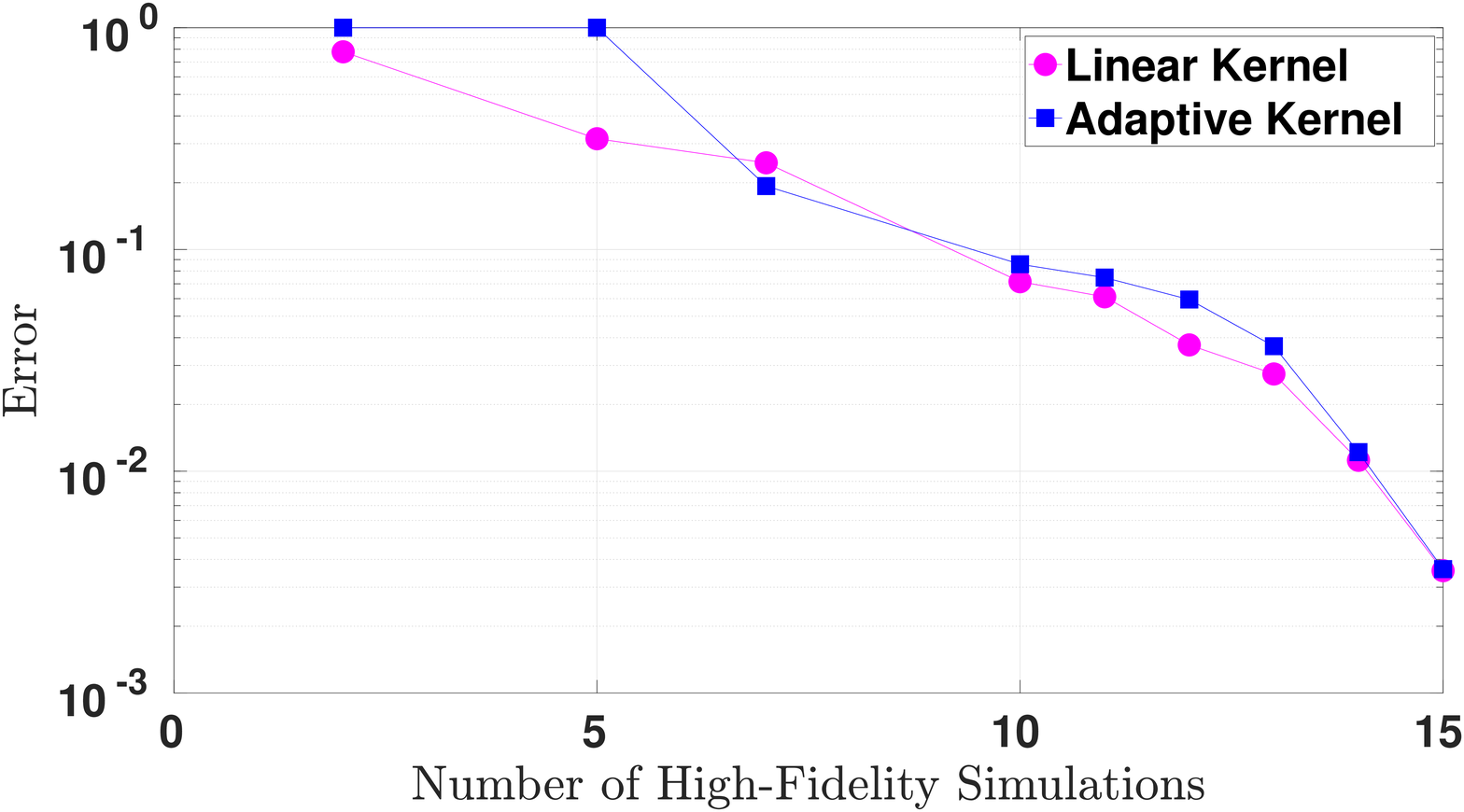}
    \caption{top wall shear stress}
  \end{subfigure}
%\subfloat[top wall shear stress]{\label{Cfig6}\epsfig{figure=figure/Stenosis_wss_t.eps,width=0.48\textwidth}} \quad
  \begin{subfigure}[t]{0.4\textwidth}
    \centering
    \includegraphics[width=\textwidth]{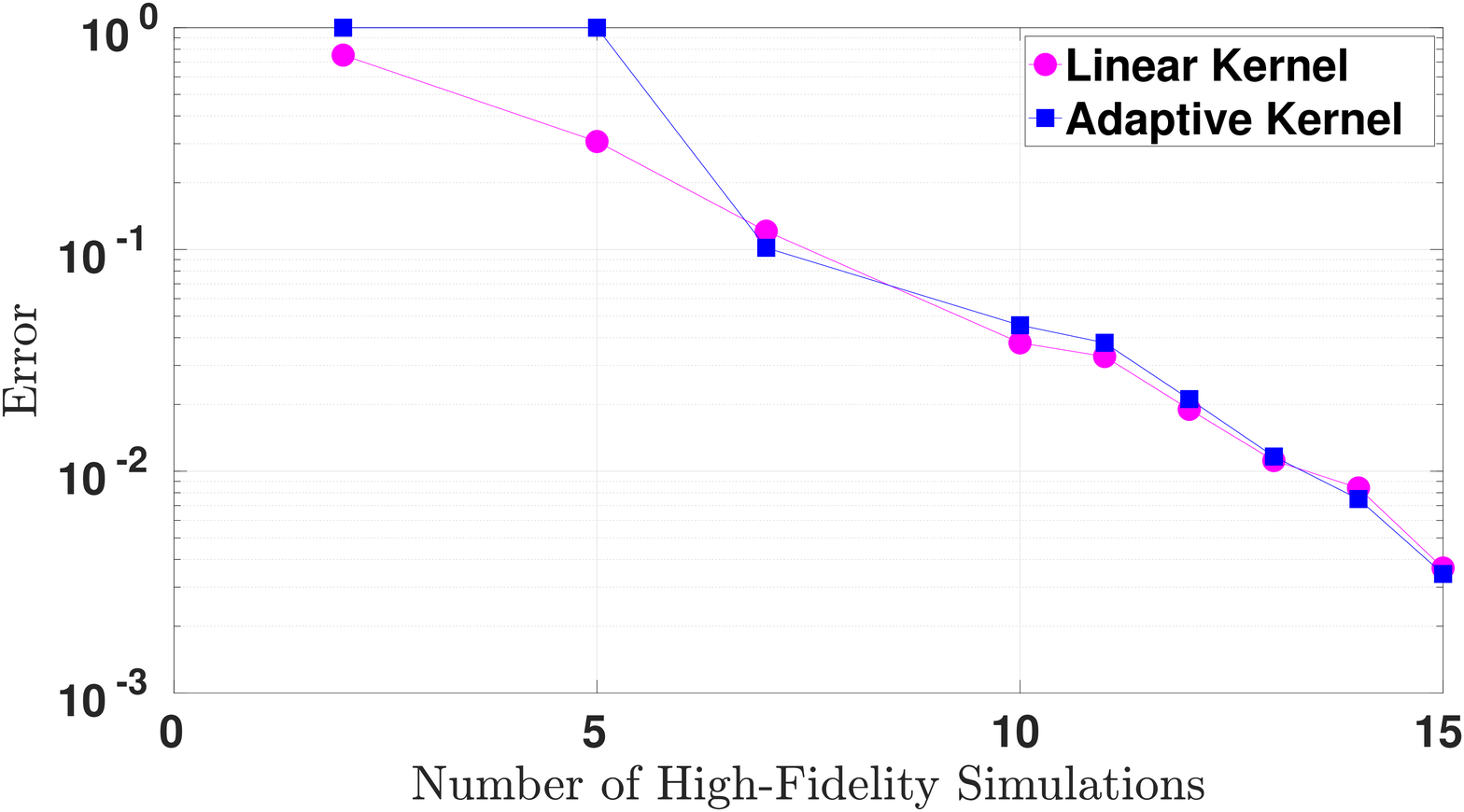}
    \caption{bottom wall shear stress}
  \end{subfigure}
%\subfloat[bottom wall shear stress]{\label{Dfig6}\epsfig{figure=figure/Stenosis_wss_b.eps,width=0.48\textwidth}} \quad
\caption{Median error of the multi-fidelity model constructed based upon the results from the full rank low-fidelity model with a coupled linearized Navier stokes equation solution in the prediction of the flow properties in a two-dimensional channel with stenosis for the corresponding high-fidelity models with full Navier-stokes equation solution; Test Problem 5; $t=30$ sec; 25 data points were tested for each case to measure the error.}
\label{fig6}
\end{figure}

\begin{figure}[!tb]
  \centering
  \begin{subfigure}[t]{0.32\textwidth}
    \centering
    \includegraphics[width=\textwidth]{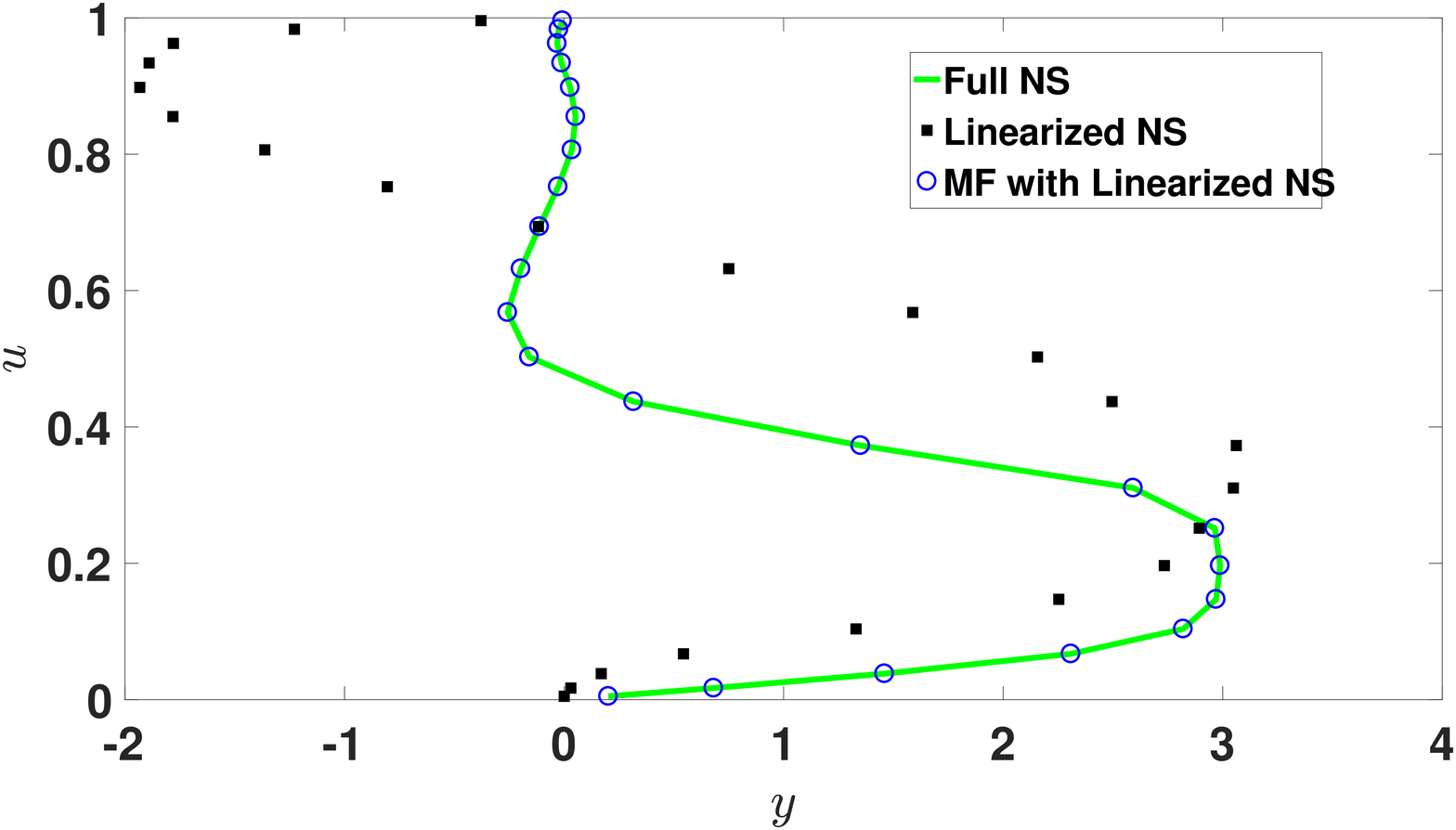}
    \caption{$u$ profile at the middle of the channel}
  \end{subfigure}
  %\subfloat[$u$ profile at the middle of the channel]{\label{Afig55}\epsfig{figure=figure/Comp_u_jj30.eps,width=0.48\textwidth}} \quad
  \begin{subfigure}[t]{0.32\textwidth}
    \centering
    \includegraphics[width=\textwidth]{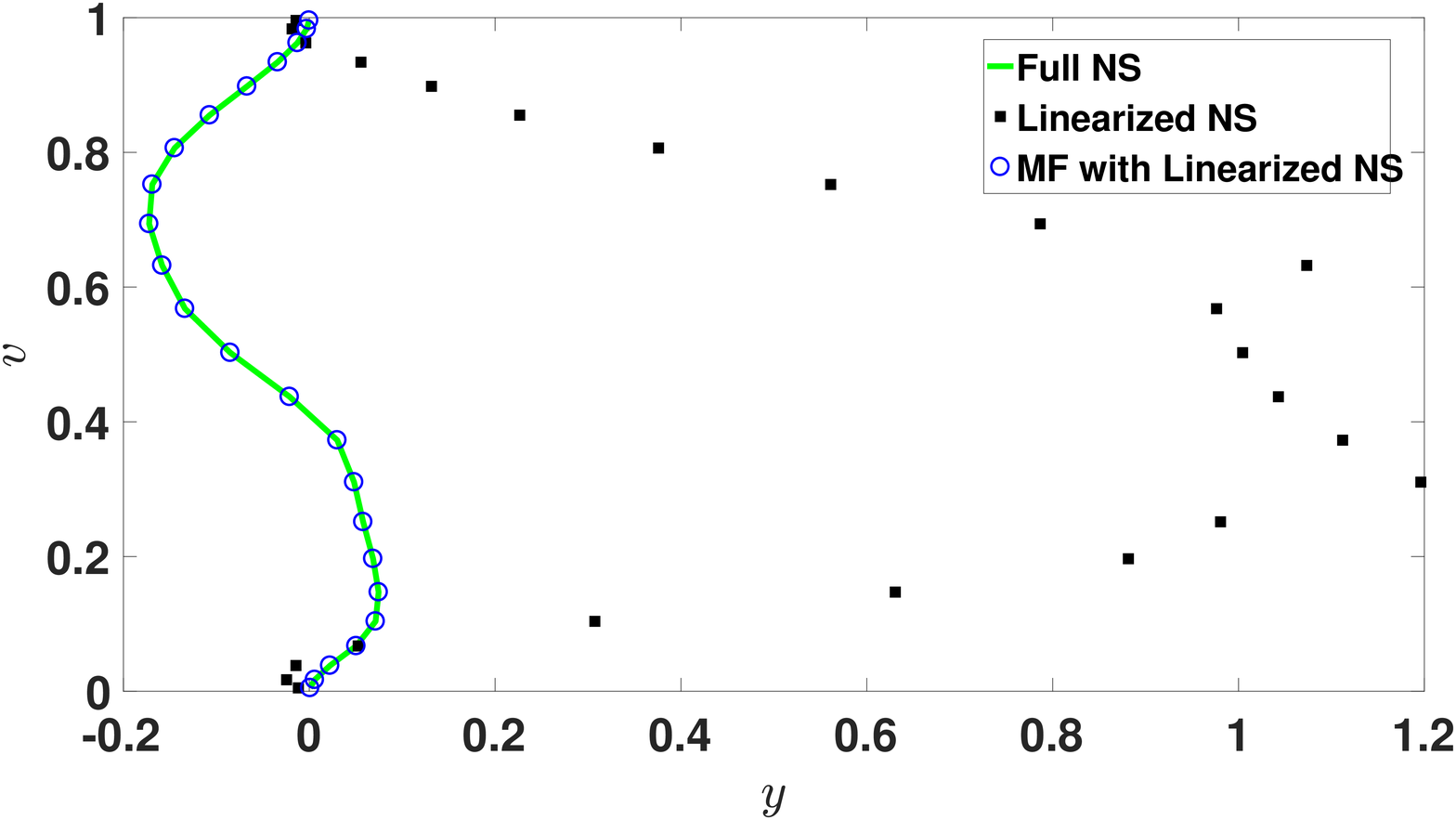}
    \caption{$v$ profile at the middle of the channel}
  \end{subfigure}
  %\subfloat[$v$ profile at the middle of the channel]{\label{Bfig55}\epsfig{figure=figure/Comp_v_jj30.eps,width=0.48\textwidth}} \newline\newline\newline\vskip -0.5cm
  \begin{subfigure}[t]{0.32\textwidth}
    \centering
    \includegraphics[width=\textwidth]{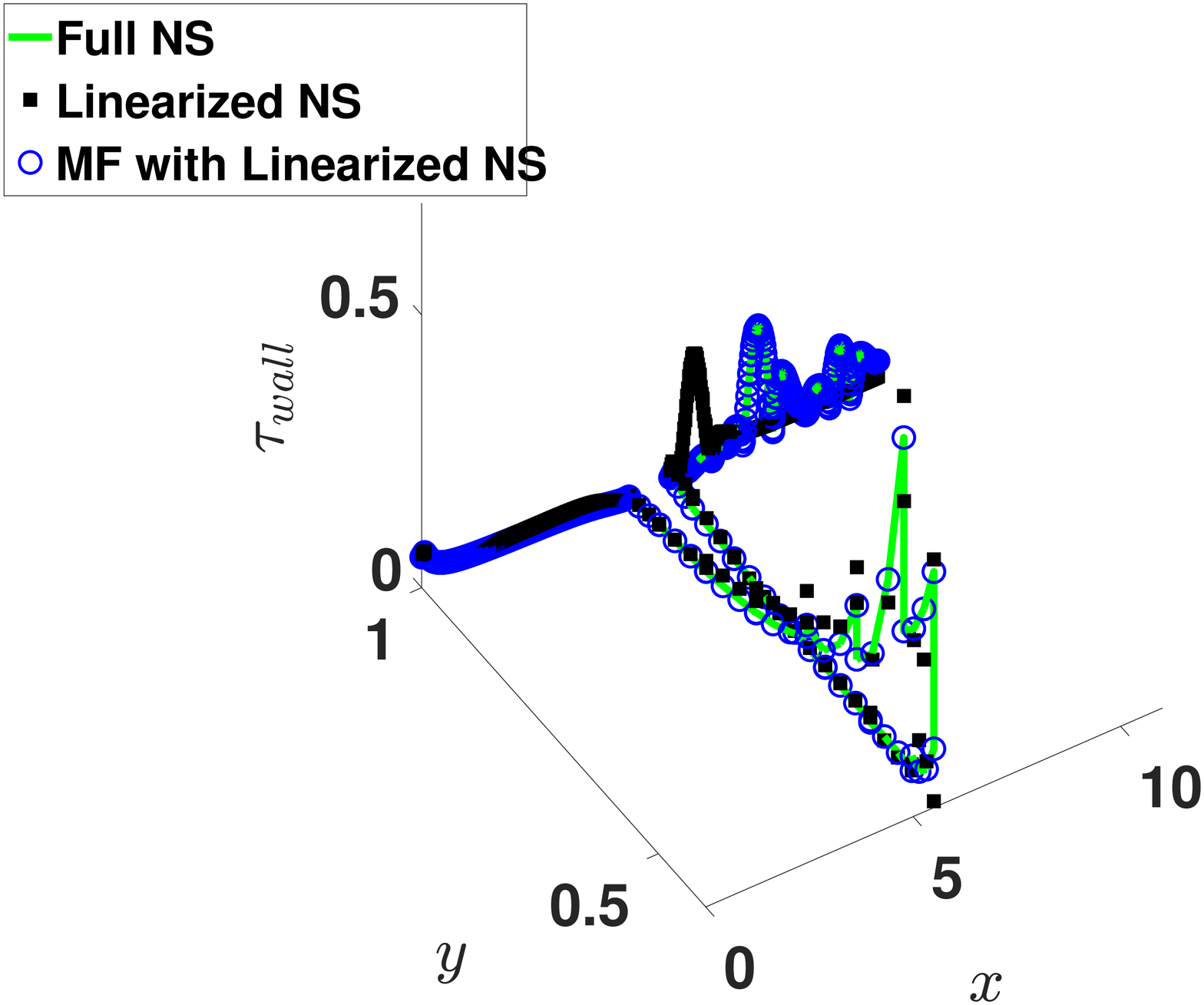}
    \caption{$\tau_{\text{wall}}$ for the top ($y=1$)}
  \end{subfigure}
  %\subfloat[$\tau_{\text{wall}}$ for the top ($y=1$)]{\label{Cfig55}\epsfig{figure=figure/Comp_twtop_jj30.eps,width=0.9\textwidth}} \quad
  \caption{Comparison of properties predicted for two-dimensional stenotic flows as obtained using high fidelity model (full Navier-stokes (NS) equation solution), low fidelity model (coupled linearized Navier stokes equation solution), and multi-fidelity approach based on low fidelity data and 16 samples from high-fidelity; Test Problem 5; $t=30$ sec, $\text{Re}=500$ and $f_0=0.6$}
  \label{fig55}
\end{figure}

Both Eqs.~\eqref{eq:P5_2} and ~\eqref{eq:P5_4} can be cast into a the form of $L(U)=B$ for a linear operator $L$ and $U$-independent term $B$. Via a Galerkin discretization, this can subsequently be cast in matrix form:
\begin{equation}
A\hat{U}=\hat{B},
\label{eq:P5_5}
\end{equation} 
where
\begin{eqnarray}
\hat{B}_m&=&\int_{\Omega}{\Phi_mBdx},\nonumber\\
A_{m,n}&=&\int_{\Omega}{\Phi_mL\left(\Phi_n\right)Bdx},\nonumber\\
\label{eq:P5_6}
\end{eqnarray}
where $\Phi_m$ are basis functions for the Galerkin method. We solve this problem using Nektar++, which is an open-source spectral/hp element software~\cite{cantwell2015nektar++}. Fourth order polynomials are used to approximate the finite element solution of the problem. Both low- and high-fidelity models use the same mesh but they are different with respect to their corresponding solution models, and as shown in Fig.~\ref{fig666} this difference leads to a significant difference in the prediction of the flow field in the channel at $t=30$ sec. The accuracy of the quantities of interest obtained from low-rank MF models with respect to the effective number of high-fidelity model evaluations is illustrated in Fig.~\ref{fig6}. In this figure, the adaptive kernel MF method produces accuracy comparable to the linear kernel MF method for all quantities of interest. Figure ~\ref{fig55} provides a clear picture for the inaccuracy of the low-fidelity model and the resulting enhancement in the predictions of the MF models. These results indicate that there is no penalty for considering a more complex kernel selection procedure in this example. Since the more complex kernel selection procedure does perform better than the linear kernel method for other examples, this suggests that using the adaptive kernel approach rather than the linear kernel approach is a better choice in general.

\section{Concluding remarks} \label{sec:conclusion}
Multi-fidelity surrogate modeling is an active area of research in the field of uncertainty quantification, with a wide range of application in different areas of science and engineering. In the context of \emph{low rank} multi-fidelity modeling, one commonly used technique is to exploit the ``kernel trick'' to express correlations in parameter space. Such a technique is the building block of low-rank MF learning algorithms and the crucial reason behind their success in providing accurate predictions of quantities of interest for complex systems. While linear kernel functions have been historically used and have demonstrated success, several drawbacks arise such as rank-deficiency and ill-conditioning of Gramian matrices when the dimension of the output space is small. The result is often numerical instability and inaccuracy of the emulator. 
   
To address this issue, two novel approaches are proposed in this paper for data-driven kernel function selection from a library existing and popular kernel functions. The selection procedure is an optimization but only relies on inexpensive low-fidelity data. Compared to the cost of a single HF evaluation, this optimization is negligible in practice. We have made some heuristic choices for hyperparameter optimization objective function \eqref{eq:Eq6} and in the choice of our optimizer (PSO with gradient descent), but alternative optimization techniques can be employed without changing the essential nature of our proposed methods. In addition, the kernel library (see Table \ref{tab:kernels}) may be grown if more kernel families are desired, or shrunk if some are deemed unnecessary. We have proposed two techniques, an Additive Kernel Approach, and an Adaptive Kernel Approach, for selection of the final kernel. 

To demonstrate the effectiveness of our methods, we have applied these new techniques to five non-trivial problems in disparate areas of science and engineering: 
molecular dynamics simulations, $n$-body galaxy models, associating polymer networks, plasmonic nano-particle arrays, and incompressible flow in two-dimensional channels with stenosis. Our numerical experiments suggest:
\begin{enumerate}
  \item the standard linear kernel is sometimes unreliable when models have a low-dimensional output space, cf. Figure \ref{fig7};
  \item \revision{the new kernel selection methods have succeeded in producing accurate and relatively inexpensive surrogates when the linear kernel fails for low-dimensional output spaces, and can achieve errors that are smaller by one-to-two orders of magnitude for identical cost, cf. Figure \ref{fig4};}
  \item the kernel selection and optimization procedure optimizes only over the low-fidelity model, and hence is very inexpensive compared to a single high-fidelity model evaluation;
  \item when the model output space is high-dimensional, the linear kernel procedure works well, but the novel kernel selection procedures do not produce worse emulators, cf. Figure \ref{fig6}; and
  \item the proposed methods can succeed on a diverse range of problems in science and engineering capability.
\end{enumerate}
Our results also suggest a small preference for the Adaptive Kernel Approach, due to a slightly better overall accuracy, along with its better computational scalability as the number $L$ of kernel functions in the library grows. For very large $L$, the complexity of a high-dimensional (convex) optimization for the alternative Additive Kernel Approach can severely impact the efficiency of this procedure. Hence, the authors generally recommend the implementation of the Adaptive Kernel Approach if the user is free to choose.

There are still open questions and challenges with this approach.  Low-rank MF methods are limited in general to cases in which the low-fidelity and high-fidelity models share a common parameter space and have similar parametric variations. Also, these procedures do not provide an approach for construction of a low-fidelity model given a trusted high-fidelity model. \revision{Finally, the procedures we have developed do not directly tackle the curse of dimensionality for high-dimensional parameter spaces.} Future work will be devoted to tackling these challenges. 

\section*{Acknowledgements}
M. Razi and R. M. Kirby acknowledge that their part of this research was sponsored by ARL under Cooperative Agreement Number W911NF-12-2-0023. The views and conclusions contained in this document are those of the authors and should not be interpreted as representing the official policies, either expressed or implied, of ARL or the U.S. Government. The U.S. Government is authorized to reproduce and distribute reprints for Government purposes notwithstanding any copyright notation herein. M. Razi and A. Narayan also acknowledge partial support from AFOSR FA9550-15-1-0467 and AFOSR FA9550-20-1-0338.

\bibliographystyle{siam}
\bibliography{bibliography}

\end{document}